\documentstyle[12pt,leqno,amscd,amssymb,verbatim]{amsart}
\newtheorem{thm}{Theorem}[section]
\newtheorem{lem}[thm]{Lemma}
\newtheorem{cor}[thm]{Corollary}

\theoremstyle{definition}

\newtheorem{defn}[thm]{Definition}

\newtheorem{rem}[thm]{Remark}
\newtheorem{rems}[thm]{Remarks}

\numberwithin{equation}{thm}

\def\sB{{\cal B}}

\def\sD{{\cal D}}\def\wsD{{{\cal D}}}
\def\sE{{\cal E}}
\def\sF{{\cal F}} 

\def\sH{{\cal H}}

\def\sO{{\cal O}}

\def\sS{{\cal S}}
\def\sT{{\cal T}}

\def\wW{{{W}}}
\def\sX{{\cal X}}
\def\sY{{\cal Y}}
\def\sZ{{\cal Z}}
                 \def\wx{{x}}
\def\dW{{\check{W}}}\def\dT{{\check{\cal T}}}\def\dS{{\check{S}}}
\def\dH{{\check{\cal H}}}\def\dx{{\check{x}}}\def\dy{{\check{y}}}
\def\dD{{\check{\cal D}}} 
\def\dC{{\check{C}}}\def\dl{{\check{\ell}}}

\def\dv{{\check{v}}}\def\ce{{\check{e}}}\def\df{{\check{f}}}
\def\tx{{\tt x}}\def\dg{{\check{g}}}
\def\dpi{{\check\pi}}

\def\bfs{{\bold s}}
\def\bft{{\bold t}}
\def\bfu{{\bold u}}

\def\fks{{\frak s}}
\def\fkt{{\frak t}}

\def\bT{{\bold T}}

\def\BQ{{\Bbb Q}}

\def\BZ{{\Bbb Z}}
\def\fT{{\frak T}}
\def\fS{{\frak S}}

\def\partial{\delta}

\def\al{\alpha} 
\def\be{\beta} 
\def\de{\delta}
\def\ep{{\varepsilon}}
\def\ga{\gamma}

\def\om{\omega} 

\def\De{\Delta}
\def\la{\lambda}\def\laa{{\lambda^{(1)}}}\def\lab{{\lambda^{(2)}}}
\def\ka{\kappa} 
\def\th{\theta} 
 
\def\La{\Lambda}

\def\oW{{\bar W}}

\font\du=cmr6
\def\dub{\text{\du\char'102}}
\def\dud{\text{\du\char'104}}
\def\oneb{{1\!\dub}}
\def\oned{{1\!\dud}}
\def\onedp{{1\!\dud'}}
\def\End{\text{\rm End}}
\def\Ind{\text{\rm Ind}} 
\def\rad{\text{\rm rad}}
 
\def\Ext{\text{\rm Ext\,}}
\def\Hom{\text{\rm Hom}}

\def\dim{\text{\rm dim\,}}

\begin{document}
\title[$q$-Schur Algebras of Type $D$]%
{Stratifying $q$-Schur Algebras of Type $D$}
\author{Jie Du}
\address{School of Mathematics, University of New South Wales,
Sydney 2052, Australia.}
\address{{\it Home page}: {\tt http://www.maths.unsw.edu.au/$\sim$jied}} 
\email{j.du@@unsw.edu.au}
\author{Leonard Scott}
\address{Department of Mathematics, University of Virginia,
Charlottesville, VA 22903, USA.}
\address{{\it Home page}: {\tt http://www.math.virginia.edu/$\sim$lls2l}} 
\email{lls2l@@virginia.edu}
\date{6 December, 1999}
\thanks{The authors
would like to thank ARC for support as well as NSF, and the Universities of 
Virginia and New South Wales for
their cooperation.} 

\begin{abstract} Two families of  $q$-Schur algebras  
associated to Hecke algebras of type $D$
are introduced, and related to a family used by Geck, Gruber and Hiss
\cite{GH}, \cite{GrH}. We  prove that the algebras
in one family, called the  $q$-Schur$^{1.5}$ algebras,
 are integrally free, stable
under base change, and  are standardly stratified if the base field
has odd characteristic. 
In the so-called linear prime case of \cite{GH},\cite{GrH},
all  three families give rise to Morita equivalent algebras.
A final section discusses a different example, and speculates
on the direction of a general theory.
\end{abstract}

\maketitle

Following up our recent work \cite{DS2} for type $B$ (and $C$), we introduce
here some endomorphism algebras of type $D$. These are of possible use in determining irreducible representations of finite groups of Lie type $D$, especially in non-defining characteristics, as in the work of Geck-Hiss \cite{GH} and Gruber-Hiss \cite{GrH}
 following
the spirit of Dipper-James' work in type $A$.
(See Remark \ref{mot} below for further discussion.)

We organize the paper as follows:
Section 1 sets up preparation on Weyl groups of classical types
and distinguished coset representatives, especially those with 
the trivial intersection property. The various $q$-Schur algebras
are introduced in Section 2 and connections between them
are also discussed. The linear prime Morita equivalence
theorems are proved in Section 3.
In Sections 4, 5 and 6, we establish further results on 
$q$-permutation and twisted $q$-permutation modules of type $B$. These results
 are supplementary
to those given in \cite{DS2}. In particular, we prove, in the type $B$
setting, the homological property
required  (see \cite{DPS1}) in stratifying an endomorphism algebra
for type $D$.
The main results are given in Section 7, where 
the quasi-heredity in the odd degree case is proved
for one of our algebras, called the $q$-Schur$^{1.5}$ algebra.
A standard stratification
in the even degree case
is constructed for the same algebra when 2 is invertible
in the base ring. Section 8 gives an effective approach to the bad prime
case $p=2$, and concludes with speculations regarding a general theory.

A preliminary version of the work was posted at the Newton Institute
workshop, June 1997, for the program on the representation theory of algebraic
and related finite groups.
\bigskip
\section{Weyl groups of classical types}
\noindent
{\bf Notation.} 
Consider the following Coxeter systems:\medskip
\halign{#\ &#\ \hfil&with\ #&#\hfil\cr
$(W_r,S)$     & the Weyl group of type $B_r$&
$S=$&$\{s_0, s_1,\cdots,s_{r-1}\}$\cr
$(\check W_r,\dS)$& the Weyl group of type $D_r$&
$\check S=$&$\{ u, s_1,\cdots,s_{r-1}\}$\cr
$(\bar W_r,\bar S)$& the Weyl group of type $A_{r-1}$&
$\bar S=$&\quad\,$\{ s_1,\cdots,s_{r-1}\}$\cr}
\medskip
We often drop the subscript $r$ and will choose
$u=s_0s_1s_0\in W$ in the sequel. Then $\oW\subseteq\dW\subseteq W$.
Let
$t_1=s_0$, $t_i=s_{i-1}t_{i-1}s_{i-1}$, $u_i=t_1t_i$ for
$2\le i\le r$, and let $C=\langle t_1,\cdots,t_r\rangle$ and
$\check C=\langle u_2,\cdots,u_r\rangle$.
Then, we have $W=C\rtimes\bar W$ and $\dW=\check C\rtimes\bar W$.
Moreover, we may identify $\oW$ with the symmetric group $\fS_r$
on $r$ letters.

Let $\ell$ (resp. $\dl$) be the length function on $\wW$ 
(resp. $\dW$) with respect to $S$ (resp. $\dS$) 
and $n_0$ the function giving
the number of $s_0$ in a reduced expression of an element of $\wW$
(see \cite[\S2.1]{DS2}). It is well-known that $\wW$ identifies with
a subgroup of the symmetric group $\fS_{2r}$ (see, e.g., \cite[(2.1W3)]{DS2}).
In this identification, the restriction to $\wW$ of the signature
function on $\fS_{2r}$ induces the group homomorphism 
$\rho_0:\wW\to\{1,-1\}$ defined by
$\rho_0(w)= (-1)^{n_0(w)}$ for all $w\in\wW$. The kernel 
of $\rho_0$, $\ker(\rho_0)$, is the intersection of $\wW$ with the alternating group
${\frak A}_{2r}$. Clearly, the subgroup $\ker(\rho_0)$ is generated by the set 
$\{s_0s_1s_0,s_1,\cdots,s_{r-1}\}$, and is isomorphic to the Weyl group
$\dW$ of type $D_r$ by the map sending $u$ to $s_0s_1s_0$ and $s_i$ to $s_i$.
This agrees with our identification above of $\dW$ with a subgroup of $\wW$. 
For convenience, we call $\dC=\langle
u_2,\cdots,u_r\rangle$ the {\it bottom} part of $\dW$, where $u_i=t_1t_i$ for
$2\le i\le r$. Note that the subgroup generated
by $s_2,\cdots,s_{r-1}$ and $\dC$ is the Weyl group of type $B_{r-1}$.
Let $f$ be the automorphism of $\dW$ induced by flipping the
Coxeter graph. Then, $f$ on $\dW$ is now the restriction
to $\dW$ of the inner automorphism $g\mapsto s_0gs_0$ of $\wW$.
Note that $f$ fixes each element of $\dC$, and interchanges
the two parabolic copies of $\fS_r$.

Fix positive integers $n,r$, and let
$\La(n,r)$ (resp. $\La^+(n,r)$) be the set of compositions
(resp. partitions) of $r$ with
$n$ parts (counting zeros in both cases).
Thus, $\La^+(r)=\La^+(r,r)$ is the set of partitions of $r$. 
Let
$$\Pi_2=\Pi_2(n,r)=\cup_{r_1+r_2=r}\La(n_r,r_1)\times\La(n,r_2),$$ 
where
$n_r$ is the maximum of $n$ and $r$.
The elements $\lambda$ in $\Pi_2$ will be
written in the sequel as $\la=(\la^{(1)},\la^{(2)})$,
 and will be called bicompositions of $r$. Let 
$\Pi_2^+$ be the subset consisting of all $(\laa,\lab)\in\Pi_2$ such that
both $\laa$ and $\lab$ are partitions. The elements in $\Pi_2^+$ are called bipartitions. Put $\Pi^+(r)=\Pi^+(r,r)$. 
We define
$$\aligned
\Pi_1&=\Pi_1(n,r)=\{\la\in\Pi_2\mid |\lab|=0\}\cr
\Pi_{\oneb}&=\Pi_{\oneb}(n,r)=\{\la\in\Pi_2\mid |\laa|=0\},\cr
\endaligned$$
where $|\la^{(i)}|$ denotes the sum of the parts of $\la^{(i)}$. 
Each $\la\in\Pi_2$ gives naturally a composition $\bar\la\in\La(n_r+n,r)$ by
juxtaposition.
For a bicomposition 
$\la=(\la^{(1)},\la^{(2)})\in\Pi_2$ with $|\laa|=a$, let 
$C_\la=\langle t_1,\cdots,t_{a}\rangle$, let $\frak S_{\bar\la}$
the Young subgroup defined by $\bar\la$, and
define
$$
W_\la=C_\la\fS_{\bar\la},\qquad
\dW_\la=W_\la\cap\dW\quad\text{ and }\quad
\bar W_\la=\frak S_{\bar\la}.
$$
In the spirit of \cite{DS2}, we shall call these subgroups {\it quasi-parabolic} subgroups.
Note that $\dW_\la=\bar W_\la$ if $a\le1$, and
$\dW_\la=\langle u_2,\cdots,u_{a}\rangle\bar W_\la$ if $a\ge2$.
Note also that $W_\la$ is generated by a unique subset $J(\la)$ of
$\tilde S=S\cup\{t_1,\cdots,t_r\}$, and hence, $J$ defines a 
map from $\cup_{n\ge 0}\Pi_2(n,r)$ to the power set $P(\tilde S)$
of $\tilde S$.

For $W_\la$, the notion of `distinguished'
coset representative\footnote{A distinguished 
coset representative $d$ is the unique element of minimal
length in $dW_\la$ or $W_\la d$. It is not necessarily minimal
with respect to the Bruhat order if $W_\la$ is not parabolic.
See \cite{Ma} or \cite{DS2} for more details.}
was introduced in \cite[\S2]{DS2}. Following \cite{DS2},
$\sD_\la$ denotes the set of distinguished right $W_\la$-coset
representatives, and $\sD_{\la,\mu}=\sD_\la\cap\sD_\mu^{-1}$
the set of distinguished double $W_\la$-$W_\mu$-coset
representatives. Let $W_{\hat\la}$ be the minimal parabolic
subgroup containing $W_\la$. (Thus, $\hat\la=(|\laa|,\lab)$.)
Then, for $d\in\sD_{\la,\mu}$,  we may write $d=u\hat dv$,
where  $\hat d\in W_{\hat\la}d W_{\hat\mu}$ is distinguished
(for parabolic subgroups), $u\in W_{\hat\la}$ and $v\in 
\sD_{W_{\hat\la}^d\cap W_{\hat\mu}}
\cap W_{\hat\mu}$. This is called
the Howlett {\it right distinguished
decomposition} of $d$ (see \cite[2.3]{DS2}).
 Putting  $\bar \sD_\la=\sD_\la\cap\bar W$ and
$\bar\sD_{\la,\mu}=\sD_{\la,\mu}\cap\bar W $, we obtain 
distinguished $\bar W_\la$-coset and double $\bar W_\la$-$\bar W_\mu$-coset
representatives. Let
$$\sD_{\la,\mu}^0=\{d\in\sD_{\la,\mu}\mid W_\la^d\cap W_\mu=\{1\}\}$$
be the set of double coset representatives in $\sD_{\la,\mu}$ 
with the {\it trivial intersection property}.
For non-negative integer $a$, let $\om_a$ denote the bipartition
$(1^a,1^{r-a})$. Then $\hat \om_a=(a,1^{r-a})$ and $W_a=W_{\hat \om_a}$.

\begin{lem}\label{tip} Let $a,b$ be non-negative integers.

$(a)$ The set $\sD_{\om_a,\om_b}^0$ is empty unless $a+b\le r$.

$(b)$ For $x\in\oW_{(a,r-a)}$, $y\in\oW_{(b,r-b)}$, we have
$d\in\sD_{\om_a,\om_b}^0$ if and only if $xdy\in\sD_{\om_a,\om_b}^0$.

$(c)$ For $d\in W$, $d\in \sD_{\om_a}$ if and only if
$n_0(xd)=n_0(x)+n_0(d)$ for all $x\in C_{\om_a}$, and
$d\in\sD_{\om_a,\om_b}^0$ if and only if $n_0(xdy)=n_0(x)+n_0(d)+n_0(y)$
for all $x\in C_{\om_a}$ and $y\in C_{\om_b}$

$(d)$ For every choice $\ep_j\in\{0,1\}$ with 
$1\le j\le i$, we have
$$t_{a+1}^{\ep_1}\cdots t_{a+i}^{\ep_i}\sD_{\om_{a+i},\om_b}^0\subseteq \sD_{\om_a,\om_b}^0.$$
\end{lem}

\begin{pf} The statement (a) is obvious, and (b) follows from 
\cite[(2.2.7)]{DS2} and the fact that $C_{\om_a}^d\cap C_{\om_b}=\{1\}
$ if and only if $C_{\om_a}^{xdy}\cap C_{\om_b}=\{1\}$
 for all $x\in\oW_{(a,r-a)}$, $y\in\oW_{(b,r-b)}$. 
We now prove (c).
Write $d=wd'$ with $w\in W_a$ and $d'\in\sD_{\hat\om_a}$.
Then $d\in \sD_{\om_a}$ if and only if $w\in\oW_a$, which is equivalent
to $n_0(xd)=n_0(x)+n_0(d)$ for all $x\in C_{\om_a}$, proving the
first part.
Now, let
$d=u\hat dv$ be the Howlett decomposition of $d$, where
 $\hat d$ is the distinguished double coset representative
in $W_a dW_b$ (for parabolic subgroups), $u\in W_a$ and $v\in W_b$. Then, 
by (b), 
$d\in\sD_{\om_a,\om_b}^0$ iff $u\in \oW_a$ and $v\in \oW_b$ (or
$d\in\sD_{\om_a,\om_b}$) and
$\hat d\in\sD_{\om_a,\om_b}^0$. The latter is equivalent to 
$s_0\not\in W_a^{\hat d}\cap W_b$,
which means no cancellation for $s_0$ when writing the 
products $x\hat dy$ ($x\in W_a$ and $y\in W_b$) in reduced form.
 (When $s_0$ is in the intersection, we have $s_0\hat d=\hat d s_0$,
and cancellation could occur.)
From this together with the first part of (c), we obtain that
$d\in\sD_{\om_a,\om_b}^0$ iff
$n_0(x dy)=n_0(x)+n_0(d)+n_0(y)$ for all $x\in C_{\om_a}$ and 
$y\in C_{\om_b}$, proving (c).
To prove (d), let $w\in\sD_{\om_{a+i},\om_b}^0$. Putting $t=t_{a+1}^{\ep_1}\cdots t_{a+i}^{\ep_i}$, 
we have, for $x\in C_{\om_a}$ and $y\in C_{\om_b}$,
$n_0(xtwy)=n_0(xt)+n_0(w)+n_0(y)$ by (c) since $w\in\sD_{\om_{a+i},\om_b}$,
But $n_0(xt)=n_0(x)+n_0(t)$, and $n_0(t)+n_0(w)=n_0(tw)$.
Therefore, $n_0(xtwy)=n_0(x)+n_0(tw)+n_0(y)$ and therefore,
$tw\in\sD_{\om_a,\om_b}^0$ by (c) again.
\end{pf}

\begin{cor}\label{tip11} Let $a\ge1$, and assume
that $d\in\sD_{\om_a,\om_b}^0$. Then 
$$\{d,t_ad,s_ad,s_at_ad\}\subset\sD_{\om_{a-1},\om_b}^0.$$
If, in addition, $d\not\in\sD_{\om_{a+1},\om_b}^0\cup
t_{a+1}\sD_{\om_{a+1},\om_b}^0$, then 
every element in $\{d,t_ad,s_ad,s_at_ad\}$ is not in 
$\cup_{\ep_i=0,1}t_a^{\ep_1}t_{a+1}^{\ep_2}\sD_{\om_{a+1},\om_b}^0$.
In other words, if we put
$U=\sD_{\om_{a+1},\om_b}^0\cup
t_{a+1}\sD_{\om_{a+1},\om_b}^0$, we have 
$$
\{1,t_a,s_a,s_at_a\}(\sD_{\om_a,\om_b}^0\backslash U)
\cap(\cup_{\ep_i=0,1}t_a^{\ep_1}t_{a+1}^{\ep_2}\sD_{\om_{a+1},\om_b}^0)=\emptyset.$$
\end{cor}

\begin{pf}
The inclusion can be proved by using Lemma \ref{tip}c and the disjointness
follows easily from the facts
$C_{\om_{a+1}}^d\cap C_{\om_b}\neq\{1\}$ and
$$U=\{x\in\sD_{\om_a,\om_b}^0\mid
 C_{\om_{a+1}}^x\cap C_{\om_b}=1\}.$$
(If $x$ is in the right hand side,
then we may write $x=tx'$ where $t\in C_{\om_{a+1}}$ and
$x'\in \sD_{\om_{a+1}}$. This is just a right coset decomposition.
Next, an argument using Lemma \ref{tip}c shows that 
$x'\in\sD_{\om_b}^{-1}$. That is, $x'\in \sD_{\om_{a+1},\om_b}$,
and even $x'\in \sD_{\om_{a+1},\om_b}^0$.
 Finally, $x\in \sD_{\om_{a}}$
forces $t=t_{a+1}$ or 1.)
\end{pf}

\section{$q$-Schur algebras of classical types}

We start with the definition of the $q$-Schur algebra of
type $B$, that is, the $q$-Schur$^2$ algebra.

Let $\sZ_0=\BZ[q,q^{-1}, q_0, q_0^{-1}]$ be 
the ring of Laurent polynomials in two
variables, and let $\sH=\sH_{q_0,q}$ be the 
(two-parameter generic) Hecke algebra 
over $\sZ_0$ associated to $W$
with defining basis $\{T_w\}_{w\in W}$. The subalgebra
$\bar \sH$ generated by $T_{s_i}$, $1\le i\le r-1$, is the Hecke algebra
associated to $\oW$. 
We denote by $\sH'=\sH_{q_0,q}'=\sH_{\sZ'}$ etc.
the Hecke algebras obtained by changing base to 
a commutative $\sZ_0$-algebra $\sZ'$. Also, for any subset $X$ of $W$,
 we denote by $\sH'(X)$
the $\sZ'$-submodule generated by all $T_w$, $w\in X$.

Define $\pi_0=\pi_0^-=1$, and for $a\ge 1$, let
$$\pi_a=\prod_{i=1}^{a}(q^{i-1}+T_{t_i}),\quad
\pi_a^-=\prod_{i=1}^{a}(q_0q^{i-1}-T_{t_i}).$$
Note that there is a $\BZ[q,q^{-1}]$-algebra automorphism on $\sH'$:
$$
\eta:\sH'\to\sH'; q_0\mapsto q_0^{-1},T_{s_0}\mapsto -q_0^{-1}T_{s_0}
\text{ and }T_{s_i}\mapsto T_{s_i}, \forall i\ge 1.
$$
Then, $\eta(\pi_a)=q_0^{-a}\pi_a^-$ for each  $a\ge0$.

The $q$-permutation $\sH'$-modules associated to $\wW_\la$, $\la\in\Pi_2$,
 are defined as cyclic
$\sH'$-modules $\wx_\la\sH'$ and $\sH'\wx_\la$ with generators 
$$\wx_\la=\pi_a \bar x_{\la},
\text{ where }a=|\laa|\text{ and } \bar x_{\la}=\sum_{w\in \oW_\la}T_w$$

\begin{defn}(\cite{DS2}) The endomorphism algebra
$$\sS_q^2(n,r;\sZ')=\End_{\sH'}(\oplus_{\la\in\Pi_2}x_\la\sH')$$
is called the $q$-Schur$^2$ algebra of degree $(n,r)$.
\end{defn}

\begin{rem}\label{djm} (a) A Morita equivalent version of the
$q$-Schur$^2$ algebra, called the $(Q,q)$-Schur algebra,
was also introduced by Dipper-James-Mathas in \cite{DJM}.

(b) Replacing $\Pi_2$ by $\Pi_1$ in the definition above,
we see from \cite[(6.3.1)]{DS2} that the centralizer 
subalgebra $\sS_q^1(n_r,r;\sZ')$ of 
the $q$-Schur$^2$ algebra 
defined by using $\Pi_1$ is isomorphic to (and hence, will be identified with)
the usual $q$-Schur algebra 
$\sS_q(n_r,r;\sZ')$.
(Recall $n_r=\text{max}(n,r)$.) This algebra was used by Dipper and James
to parametrize the irreducible modular characters
of finite GL$_n$ in the non-defining characteristic case,
while the one defined by using $\Pi_{\oneb}$ was considered in
the work of Geck-Hiss and Gruber-Hiss
on finite orthogonal and symplectic groups. 
We call this latter algebra here the $q$-Schur$^{\oneb}$
algebra.
Note that its identity is an idempotent in the $q$-Schur$^2$ algebra.
We agree to not insist on the same identity element when using
the terminology `centralizer subalgebra'.
Also, we shall write $B\le A$ if $B$ is a centralizer subalgebra of $A$.
That is, $B=eAe$ where $e$ is an idempotent in $A$.
Recall from \cite[(6.3.2)]{DS2} that, if  $B\le A$ where $A$ is an
$\sO$-free algebra for a regular ring $\sO$, 
then the decomposition matrix of $B$ is part of the decomposition matrix of 
$A$. 

(c) Using the twisted permutation module $y_\la\sH'$, where
$$
y_\la=\pi_{|\laa|}^-\bar y_\la,\text{ with }
\sum_{w\in\bar W_\la}(-q)^{-\ell(w)}T_w,
$$
 we may define another endomorphism algebra $\tilde\sS_q^2(n,r;\sZ')$.
It is known \cite{Du1} that $\tilde\sS_q^2(n,r;\sZ')$ is isomorphic
to $\sS_q^2(n,r;\sZ')$ as $\sZ'$-algebras.
\end{rem}

\begin{thm}\label{qh} (\cite{DS2}) The $q$-Schur$^2$ algebra
$\sS_q^2(n,r;\sZ')$ is (integrally) 
quasi-hereditary.
\end{thm}

In the so-called linear prime case (a fairly strong restriction),
the following result was first obtained by Gruber and Hiss \cite{GrH}
in the case $n=r$.
It follows from (\ref{qh}) and \cite[(6.3.2)]{DS2} in general.

\begin{cor} 
The decomposition matrix for the $q$-Schur$^{\oneb}$
algebra contains an $m\times m$ upper unitriangular block, where
$m$ is the number of non-isomorphic irreducible modular representations.
\end{cor}

In the linear prime case studied in  \cite{GH} and \cite{GrH}, the decomposition
matrix is square. That is, $m$ is the same as the number of 
ordinary irreducible
representations (the number of bipartitions).
The validity of the conjecture \cite[3.4]{GH} would imply this is true
for all odd primes.\footnote{In the context of \cite{GH},\cite{GrH},
$q$ is a power of a rational integer prime, which is reduced modulo
a different prime when discussing decomposition numbers.}

We now turn to defining the $q$-Schur$^{1.5}$, $q$-Schur$^{2.5}$ and
$q$-Schur$^{\oned}$ algebras
associated to the Hecke algebra of type $D$.

Let $\dH$ be the (generic) Hecke algebra over 
$\sZ=\BZ[q,q^{-1}]$ associated to the Weyl group $\dW$,
and let $\dH'=\dH\otimes_{\sZ}\sZ'$ for any commutative
$\sZ$-algebra $\sZ'$. 
We define (compare \cite{P})
$$\dpi_0=\dpi_1=1,\quad \dpi_a=\prod_{i=2}^a(q^{i-1}+T_{u_i}), \quad 
2\le a\le r,$$
and put $\dx_\la=\dpi_a \bar x_{\la}$, where $a=|\laa|$.
Then the $\dH'$-modules $\dx_\la\dH'$ and $\dH' \dx_\la$ are called
the {\it $q$-permutation modules} associated to the subgroup $\dW_\la$.

For $\la\in\Pi_{\oneb}$ (i.e. $|\laa|=0$), $f(\dW_\la)$ is a
parabolic subgroup, and is different from $\dW_\la$ (not even conjugate to it)
if $s_1\in\dW_\la$.  
We define the notation $f(\la)$ by $\dW_{f(\la)}=\dW_{f(J(\la))}=f(\dW_\la)$ and put
$$\aligned
\Pi_{2.5}&=\Pi_{2.5}(n,r)=\Pi_2(n,r)\cup f\Pi_{\oneb}(n,r)\cr
\Pi_{1.5}&=\Pi_{1.5}(n,r)=\{\la\in\Pi_2(n,r)\mid |\laa|\ge|\lab|\}.\cr
\endaligned
$$

\begin{defn} \label{typeD} For any positive integers $n,r$ with $r\ge4$, let
 $\dT^\ka_{\sZ'}=
\oplus_{\la\in\Pi_\ka}\dx_\la\dH'$ for $\ka=1.5,2.5$. 
The endomorphism algebra
$$\sS_q^\ka(n,r;\sZ')=\End_{\dH'}(\dT^\ka_{\sZ'})$$
is called the {\it $q$-Schur$^\ka$ algebra} of degree ($n,r$).
By using the disjoint union $\Pi_{\oned}=\Pi_{\oneb}\sqcup f(\Pi_{\oneb})$,
we similarly define $\dT_{\sZ'}^{\oned}$.  The corresponding centralizer
subalgebra of the $q$-Schur$^{2.5}$ algebra,
an analog for type $D$ to the $q$-Schur$^{\oneb}$ algebra
for type $B$,
is called a {\it $q$-Schur$^{\oned}$ algebra}.
When $\sZ'=\sZ$, we will denote $\sS_q^\ka(n,r;\sZ)$ simply by
$\sS_q^\ka(n,r)$.
\end{defn}

\begin{rems} \label{onedp} 
(1) Using (\ref{cli}b) below, we see that the $q$-Schur$^{\oned}$ 
algebra is isomorphic to the endomorphism algebra
$$\End_{\dH'}(\oplus_{\la\in \Pi_\oneb}\dx_\la\sH'|_{\dH'}).$$

(2) The $q$-Schur algebra of type $D$ defined
in \cite[7.2]{GrH} uses only one parabolic copy of $\fS_r$.
In terms of our notation, it is (Morita equivalent to)
the endomorphism algebra
of the $\dH'$-module $\oplus_{\la\in\Pi_\oneb}\dx_\la\dH'$
(or $\oplus_{\la\in f\Pi_\oneb}\dx_\la\dH'$).
We shall denote this algebra by $\sS_q^\onedp(n,r;\sZ')$.
Clearly, $\sS_q^\onedp(n,r;\sZ')\le \sS_q^\oned(n,r;\sZ')$
and therefore, its decomposition matrix is determined
by that of the $q$-Schur$^{\oned}$ algebra.
Actually, Gruber and Hiss,  later in their paper, appear
to be implicitly using $\sS_q^{\oned}(r,r)$ rather than
$\sS_q^\onedp(r,r)$. See Remark \ref{onedd}(2) below.
\end{rems}

When the $\sZ_0$-algebra
$\sZ'$ has the property that $q_0$ is specialized to 1, that is,
$\sH'=\sH_{1,q}'=\sH\otimes_{\sZ_0}\sZ'$, we have in $\sH'_{1,q}$
$$
T_{s_0}^2=1 \text{ and } T_{s_0}T_w=T_{s_0w},\quad \forall w\in \dW.
$$

{\it In the rest of the section, we assume} $\sH'=\sH'_{1,q}$. 
The first part of the following has been observed by Gruber-Hiss.

\begin{lem}\label{cli} (a) The algebra $\dH'$ is isomorphic to the subalgebra
of $\sH'=\sH'_{1,q}$ generated by $T_s$, 
$s\in\dS$. So we identify $\dH'$ with this
subalgebra and obtain a graded Clifford system for $\sH'$ over
$\dH'$.

(b) We have $\dH'$-module decompositions
$$\wx_\la\sH'=\cases x_\la\dH',&\text{ if }|\laa|\neq0\cr
             x_\la\dH'\oplus T_{s_0}\check f(x_\la)\dH',&\text{ if }|\laa|=0,\cr
\endcases$$
where $\check f$ is the induced flipping automorphism on $\dH'$ defined by
$\check f(T_w)=T_{s_0}T_wT_{s_0}=T_{f(w)}=T_{\check f(w)}$.

(c) For any $1\le a\le r$, we have 
$$\pi_a=(1+T_{s_0})\dpi_a=\dpi_a(1+T_{s_0}),\qquad
\pi_a^-=(1-T_{s_0})\dpi_a=\dpi_a(1-T_{s_0}).$$
Thus, $\dx_\la\dH'\cong\wx_\la\dH'=(1+T_{s_0})
\dx_\la\dH'$,  for any $\la\in\Pi_2$ with $|\laa|\ge1$. 
In particular, the permutation modules $\dx_\la\dH'$ are $\sZ'$-free
and pure in $\dH'$.
\end{lem}

\begin{pf} The statements (a) and (b) are clear, noting $\sH'=\dH'\oplus
T_{s_0}\dH'$.
The first assertion in (c) follows from induction and the relation
$$(1+T_{s_0})(q^{i-1}+T_{t_i})=(q^{i-1}+T_{t_i})(1+T_{s_0})=(1+T_{s_0})(q^{i-1}+T_{u_i})\quad\forall i,$$
noting $u_i=s_0t_i=t_is_0$  and $T_{t_i}=T_{s_0}T_{u_i}$ for all $i\ge2$.
For the second displayed equation in (c), use
$$(1-T_{s_0})(q^{i-1}-T_{t_i})=(1-T_{s_0})(q^{i-1}+T_{u_i}),i\ge2.$$
The freeness follows from (b) and \cite[(3.2.2a)]{DS2}, and the purity
follows from the fact that $x_\la\dH'=x_\la\sH'$ is pure in $\sH'
=(1+T_{s_0})\dH'\oplus T_{s_0}\dH'$, and hence pure in $(1+T_{s_0})\dH'$.
\end{pf}

Using the decomposition for $x_\la\sH'$ above and noting
$\sS_q^1(n,r;\sZ')\le\sS_q^1(n_r,r;\sZ')$, we obtain immediately
the following.

\begin{cor} \label{clii} Let $\sT^2_{\sZ'}=\oplus_{\la\in\Pi_2}x_\la\sH'$. Then
$\sT^2_{\sZ'}|_{\dH'}\cong\dT^{2.5}_{\sZ'}$. Hence,
the $q$-Schur$^{2.5}$ algebra is isomorphic to 
$\End_{\dH'}(\sT^2_{\sZ'})$, and  we have 
$$\cases(1)& \sS_q^1(n,r;\sZ')\le
           \sS^{1.5}_q(n,r;\sZ') \le\sS^{2.5}_q(n,r;\sZ')\cr
(2)&\sS_q^1(n,r;\sZ')\le
           \sS^{2}_q(n,r;\sZ') \subseteq\sS^{2.5}_q(n,r;\sZ')\cr
(3)& \sS^{\oned}_q(n,r;\sZ') \le\sS^{2.5}_q(n,r;\sZ').\cr
\endcases$$
\end{cor}

In this paper, we will mainly investigate  the structure
of the $q$-Schur$^{1.5}$ algebras, though we keep the $q$-Schur$^{\oned}$ algebras in mind. See
(\ref{mot}), (\ref{mrii}), (\ref{char2}b) and (\ref{badp}). 

Specht modules $S_{\la K}$ over the quotient field $K$ of $\sZ$ is
a complete set of irreducible $\sH_K$-modules. Thus, their restrictions to 
$\dH_K$ are $\dH_K$-modules. The following result is an easy
consequence of Tits' deformation theory and standard facts about ordinary
characters of Weyl groups of type $D$.

\begin{lem} \label{6a}
Let $r$ be odd. Then the set $\{S_{\la K}\mid
\la\in\Pi^+_{1.5}(r,r)\}$ is a complete set of distinct irreducible
$\dH_K$-modules. If $r$ is even, then $S_{\la K}$,
 for a bipartition $\la=(\al,\be)$ of $r$, is irreducible iff
$\al\neq\be$, and $S_{(\al,\al)K}$ has two (distinct) composition factors.
Moreover, we have in general $S_{\la K}\cong
S_{\la^\star K}$ on $\dH'$ for $\la^\star=(\be,\al)$.
\end{lem}

Recall the automorphism $\eta$ defined at the beginning of this section,
and note that $\eta(T_{s_0})=-T_{s_0}$. (Also, $\eta(T_{u_i})=T_{u_i}$
for all $i\ge 2$, and the automorhism $\eta$ fixes $\dH'$.)
Let $S_\la^\eta$ be the module $S_\la$ with $\sH$-action twisted by $\eta$.
We abbreviate below the induction functor $\Ind_{\dH'}^{\sH'}$ to
$\uparrow^{\sH'}$.

\begin{lem} \label{6b}
(a) For any $\la\in\Pi^+(r)$, let $\la^\star=(\lab,\laa)$.
Then we have $S_{\la K}^\eta\cong S_{\la^\star K}$ as $\sH_K$-modules.
Moreover, if $|\laa|\ge|\lab|$ and  $|\mu^{(1)}|\ge
|\mu^{(2)}|$, but not both equalities,
then $\Hom_{\sH_K}(x_\la^\eta\sH_{ K},S_{\mu K})=0$;
in particular, when $r$ is odd, we always have
$\Hom_{\sH_K}(x_\la^\eta\sH_{ K},x_\mu \sH_K)=0$.

(b) If $|\laa|=0,$ then $\dx_\la\dH'\uparrow^{\sH'}\cong x_\la\sH'$.

(c) If $2$ is invertible in $\sZ'$, then
we have $\dx_\la\dH'\uparrow^{\sH'}\cong x_\la\sH'\oplus x_\la^\eta\sH'$
for any $\la\in\Pi_2$ with $|\laa|\ge1$.
\end{lem}

\begin{pf}  
By the discussion in \cite[(5.1.2)]{DS2}, it suffices to look at the
specialization at $q_0=q=1$. In this case, the construction given at the 
beginning of \cite[\S5.1]{DS2} gives the required isomorphism.
Using this, the composition factors $S_{\nu K}$ appearing in $x_\la^\eta\sH_K
\cong(x_\la\sH_K)^\eta$ (the module $x_\la\sH_K$
with action twisted by $\eta$)  
have the property $|\nu^{(1)}|\le|\nu^{(2)}| $. 
So the last assertion in (a) follows.
The statement (b) follows from the fact $\dx_\la=\bar x_\la=x_\la$
in this case.
Finally, since 2 is invertible, we have $\sH'=(1+T_{s_0})\sH'\oplus (1-T_{s_0})\sH'$, and therefore,
$$\dx_\la\dH'\uparrow^{\sH'}=\dx_\la\otimes_{\sH'}(1+T_{s_0})\sH'
\oplus \dx_\la\otimes_{\sH'}(1-T_{s_0})\sH'\cong x_\la\sH'\oplus x_\la^\eta\sH'.$$
\end{pf}

\begin{rem}\label{mot}
(1) Similar to (\ref{djm}c), we define for a bicomposition $\la$
$\dy_\la=\dpi_{|\laa|}\bar y_\la$. Then, we may establish a parallel
theory as above with all $\dx_\la$ replaced 
by $\dy_\la$. In particular, the resulting $q$-Schur algebras
are isomorphic copies of the original ones and all $\dy_\la\dH'$
are pure in $\dH'$.

(2) 
Let $G$ denote the finite orthogonal group SO$_{2n}(q)$ of 
even degree and fix a splitting $p$-modular system $\{K,\sO,k\}$.
Then $\dH_{\sZ'}$ ($\sZ'\in\{K,\sO,k\}$)  is isomorphic to the endomorphism
algebra of the induced module $M_{\sZ'}$ of the trivial $B$-module $\sZ'$,
where $B$ is a fixed Borel subgroup. Then the purity
of the submodules $\dy_\la\dH_\sO$ of the symmetric algebra
$\dH_\sO$ over $\sO$ gives rise to  isomorphisms
by \cite[Thm1]{C}
$$\sS_q^\kappa(n,r;\sO)\cong \End_{\sO G}(\oplus_{\la\in\Pi_{\kappa}}\sqrt{\dy_\la M_{\sO}}),\text{ for } \kappa=1.5,2.5, \oned.$$

The arguments in \cite[\S 4]{GH} shows that, for $\kappa=\oned$,
 the decomposition matrix
of the $\sO G$-module $\oplus_{\la\in\Pi_{\kappa}}\sqrt{\dy_\la M_{\sZ'}}$
is part of the decomposition matrix of $\sO G$,
and is even decisive in determining the unipotent part of the latter
decomposition matrix in good cases. 
Probably the determination of this
decomposition matrix is about as hard as determining
that of the larger algebra $\sS_q^{2.5}(n,r;\sZ')$, though
not all irreducible modular representations of the latter may be required
in every case.
At the same time, it may be hoped that the decomposition matrix for
$\sS_q^{1.5}(n,r;\sZ')$ will be helpful in determining decomposition
numbers and modular irreducible representations of $\sS_q^{2.5}(n,r;\sZ')$,
since the two algebras have the {\it same}
ordinary irreducible representations. This philosophy
works well in at least two quite different cases, presented in section 3
below and in section 8.
 \end{rem}  

\section{The linear prime case: Morita equivalence theorems}

In this section, we shall prove that, in the linear prime case
(where much is already known, cf. \cite{GH}),
the $q$-Schur$^2$ and $q$-Schur$^{\oneb}$ algebras are Morita equivalent.
In the type $D$ linear prime case, both the 
 $q$-Schur$^{2.5}$ and $q$-Schur$^{\oned}$ algebras 
are Morita equivalent to the $q$-Schur$^{1.5}$ algebra.

We need some preparation at the group algebra level. For $0\le a\le r$, 
let 
$$\aligned
\ep_{a,r-a}&=2^{-r}(1+t_1)\cdots(1+t_a)(1-t_{a+1})\cdots(1-t_r)\cr
&=2^{-r}(1+t_1)\cdots(1+t_a)w_{a,r-a}(1-t_{1})\cdots(1-t_{r-a})w_{a,r-a}^{-1},
\endaligned$$
where   
$w_{a,b}=s_a\cdots s_1s_{a+1}\cdots s_2\cdots s_{a+b-1}\cdots s_b
=\left(\smallmatrix 1&\cdots&a&a+1&\cdots&a+b\cr
             b+1&\cdots &b+a&1&\cdots &b\cr
\endsmallmatrix\right).$
Note that $\oW_{(a,r-a)}w_{a,r-a}=w_{a,r-a}\oW_{(r-a,a)}$.
For a bicomposition $\la\in\Pi_2(n,r)$, we shall consider
 a decomposition for the second part of $\la$:
$$\lab=\be+\ga=(\be_1+\ga_1,\cdots,\be_n+\ga_n)$$
where $\be=(\be_1,\cdots,\be_n)\in\La(n,|\be|)$ and
$\ga=(\ga_1,\cdots,\ga_n)\in\La(n,|\ga|)$. 
Associated to such a decomposition $\lab=\be+\ga$, we have two compositions:
$$(\be|\ga)=(\be_1,\ga_1,\cdots,\be_n,\ga_n),\quad
(\be,\ga)=(\be_1,\cdots,\be_n,\ga_1,\cdots,\ga_n).$$
Write $a=|\laa|$, $b=|\be|$ and $c=|\ga|$, and 
let $w_{\be,\ga}^{(2)}\in \oW_{(1^a,r-a)}=\fS_{\{r-r_2+1,\cdots,r\}}$ 
be the (unique in $\oW_{(1^a,r-a)}$) distinguished
double coset representative\footnote{See, for example, \cite[3.4]{Du}
for a construction of such a representative.}
 for the subgroups
$\oW_{(1^a,\lab)}$ and $\oW_{(1^a,b,c)}$ of
$\oW_{(1^a,r-a)}$
such that 
$$\oW_{(1^a,\lab)}\cap
w^{(2)}_{\be,\ga}\oW_{(1^a,b,c)}(w^{(2)}_{\be,\ga})^{-1}=
\oW_{(1^a,(\be|\ga))},$$
and
$$
(w^{(2)}_{\be,\ga})^{-1}\oW_{(1^a,\lab)}w^{(2)}_{\be,\ga}\cap
\oW_{(1^a,b,c)}=
\oW_{(1^a,\be,\ga)}.$$

\begin{lem}\label{lp1}  Maintain the notation introduced above.
Let $R$ be a commutative ring containing $2^{-1}$.
Then, for any bicomposition $\la$ of $r$ with $|\laa|=a$, 
$$\tx_\la RW=\bigoplus_{\Sb\be\in\La(n,b),\ga\in\La(n,c)\cr
\be+\ga=\lab\endSb}\bar\tx_\la w_{\be,\ga}^{(2)}\ep_{a+b,c} RW$$
where $\tx_\la=\sum_{w\in W_\la}w$ and $\bar\tx_\la=\sum_{w\in\oW_\la}w$
\end{lem}

\begin{pf} Since 2 is invertible in $R$, we have all $\ep_{a,r-a}\in RC$ and
$RC=\oplus_{a=0}^r\oplus_{w\in\oW\backslash \oW_{(a,r-a)}} Rw\ep_{a,r-a}w^{-1}$.
Thus, 
$$\aligned
\tx_\la RW&\cong \Ind_{R\oW_\la C}^{RW}(\Ind_{RW_\la}^{R\oW_\la C}(
R\tx_\la))\cr
&\cong \bigoplus_{\be+\ga=\lab}
R\bar\tx_\la w^{(2)}_{\be,\ga}\ep_{a+|\be|,|\ga|}(w^{(2)}_{\be,\ga})^{-1}
\otimes_{R\oW_{(\laa,(\be|\ga))}C}RW.\cr
\endaligned$$
We certainly have
$$R\bar\tx_\la w^{(2)}_{\be,\ga}\ep_{a+|\be|,|\ga|}(w^{(2)}_{\be,\ga})^{-1}
\otimes_{R\oW_{(\laa,(\be|\ga))}C}RW\cong 
\bar\tx_\la w_{\be,\ga}^{(2)}\ep_{a+|\be|,|\ga|} RW,$$
proving the lemma.
\end{pf}

\noindent
{\bf The type $B$ case.}
Let $\sH=\sH_{q_0,q}$ be the Hecke algebra of type $B$ over $\sZ_0$ defined in \S2.
Recall from \cite[(3.8),(4.4)]{DJ} the elements $v_{a,b}=\pi_a T_{w_{a,b}}\pi_b^-$
and the polynomial 
$$g_r=g_r(q_0,q)=\prod_{i=-(r-1)}^{r-1}(q_0+q^i).$$ 
Let $\sZ_{0,g_r}$ be the ring obtained by localizing $\sZ_0$ at $g_r$.
Thus, $g_r$ is invertible in $\sZ_{0,g_r}$.
Let $\sH_{g_r}=\sH\otimes_{\sZ_0}
\sZ_{0,g_r}$. 
Then, by \cite[(3.27-8)]{DJ}, there are orthogonal idempotents
$\{e_{a,r-a}\mid 0\le i\le r\}$ in $\sH_{g_r}$ such that
$$v_{a,r-a}\sH_{g_r}=e_{a,r-a}\sH_{g_r},\text{ and }
\sH(\fS_{(a,b)})e_{a,b}=e_{a,b}\sH(\fS_{(a,b)}).$$
Also, $e_{a,r-a}$ specializes to $\ep_{a,r-a}$ in $\BQ W$.
The following lemma gives a new basis for $\sH_{g_r}$.

\begin{lem}\label{newbs}
The set 
$$\sB=\{T_xe_{a,r-a}T_y\mid 0\le a\le r,x\in\bar\sD_{(a,r-a)}^{-1},y\in\oW\}$$
forms a basis for $\sH_{g_r}$.
\end{lem}

\begin{pf}
Let $H$ be the $\sZ_{0,g_r}$-submodule generated by the set $\sB$.
Specializing both $q_0$ and $q$ to 1, we obtain an algebra homomorphism
$\sZ_{0,g_r}\to \BQ$. Since 
$\sH_{g_r}\otimes_{\sZ_{0,g_r}}\BQ\cong\BQ W$
and the image of the set $\sB$ becomes the basis 
$$\{(w\ep_{a,r-a}w^{-1})wy\mid 0\le a\le r,x\in \oW/\oW_{(a,r-a)},y\in\oW\}.$$
So $\sB$ is linearly independent.
To check that it spans $\sH_{g_r}$, it suffices to prove the
corresponding spanning property over any field $k$ which is a $\BZ_{0,
g_r}$-algebra.
Let $e_{a,r-a,k}$ be the image of $e_{a,r-a}$ in $\sH_{g_r,k}$.
Note that $\sH_{g_r,k}e_{a,r-a,k}\sH_{g_r,k}$ is spanned by
$T_xe_{a,r-a,k}T_y$ ($x\in\bar\sD_{(a,r-a)}^{-1},y\in\oW$). So,
if $\sH_{g_r,k}=\sum_{a=0}^r\sH_{g_r,k}e_{a,r-a,k}\sH_{g_r,k}=\sH_{g_r,k}e_k\sH_{g_r,k}$, where $e_k=\sum_{a=0}^re_{a,r-a,k}$,
then $\sB_k$ must span $\sH_{g_r,k}$. 

It remains to prove the equality $\sH_{g_r,k}=\sH_{g_r,k}e_k\sH_{g_r,k}$.
Suppose this is not the case. Then there is a maximal left ideal
in $\sH_{g_r,k}$ such that 
$$\sH_{g_r,k}e_k\sH_{g_r,k}\subseteq I.$$
Since the functor $M\mapsto e_kM$ defines an Morita equivalence
from the category $\sH_{g_r,k}$-{\bf mod} to the catgory
$e_k\sH_{g_r,k}e_k$-{\bf mod}, it must be that $e_k(\sH_{g_r,k}/I)\neq0$.
However, on the other hand,
$$e_k\sH_{g_r,k}=e_k\sH_{g_r,k}e_k\sH_{g_r,k}\subseteq e_kI
\subseteq e_k\sH_{g_r,k},$$ which implies that $e_k(\sH_{g_r,k}/I)=0$,
a contradiction.
\end{pf} 

We now have the following.

\begin{lem}\label{lp2}
Let $\la=(\laa,\lab)\in\Pi_2(n,r)$ be a bicomposition with $|\laa|=a$. Then
the set, denoted by $\sB_\la$,
$$
\{\bar x_\la T_{w^{(2)}_{\be,\ga}}e_{a+|\be|,|\ga|}T_w\mid
\be\in\La(n,b),\ga\in\La(n,c), \be+\ga=\lab, w\in\bar\sD_{(\laa,\be,\ga)}\},
$$
where $\bar\sD_{(\laa,\be,\ga)}$ is the set of distinguished right
$\oW_{(\laa,\be,\ga)}$-coset representatives in $\oW$, 
forms a basis for $x_\la\sH_{g_r}$.
\end{lem}

\begin{pf} Write 
$\bar x_\la=\bar x_\laa h_{\be,\ga}x_{\be,\ga}'$
with $x_{\be,\ga}'=\sum_{w\in\oW_{(1^a,(\be|\ga))}}T_w$
and $h_{\be,\ga}=\sum_{d\in\oW_{(1^a,\lab)}\cap\bar\sD_{(1^a,(\be|\ga))}^{-1}}T_d$.
Since $w_{\al,\be}^{(2)}$ is distinguished, we have
$\bar x_\la T_{w_{\be,\ga}^{(2)}}=h_{\be,\ga}T_{w_{\be,\ga}^{(2)}}
\bar x_{(\laa,\be,\ga)}$,
and, for each $d\in\oW_{(1^a,\lab)}\cap\bar\sD_{(1^a,(\be|\ga))}^{-1}$,
$T_dT_{w_{\be,\ga}^{(2)}}=T_{dw_{\be,\ga}^{(2)}}$ and
$dw_{\be,\ga}^{(2)}\in\bar\sD_{(a+b,c)}^{-1}$.
Note also that the union of all $(\oW_{(1^a,\lab)}\cap\bar
\sD_{(1^a,(\be|\ga))}^{-1})w_{\be,\ga}^{(2)}$ with $\lab=\be+\ga$ is disjoint.
 Thus, for any field which is a $\sZ_{0,g_r}$-algebra, the image of
$\sB_\la$ in $x_\la\sH_{g_r,k}$ is linearly independent
by (\ref{newbs}).
Therefore, if  $M$ denotes the $\sZ_{0,g_r}$-submodule generated by the set $\sB_\la$, then the image $M_k$ of $M$ in $x_\la\sH_{g_r,k}$
has dimension $|W/W_\la|$ by (\ref{lp1}),
which is equal to dim\,$x_\la\sH_{g_rk}$ by 
\cite[(3.2.2a)]{DS2}. Therefore, $M_k=\sH_{g_rk}$ for any field $k$
which is a $\sZ_{0,g_r}$-algebra. Consequently
$M=x_\la\sH_{g_r}$ and hence $\sB_\la$ forms a basis for $x_\la\sH_{g_r}$.
\end{pf}

\begin{cor} For any $\la\in\Pi_2(n,r)$ with $|\laa|=a$, we have
the following decomposition
\begin{equation}\label{lpp}
x_\la\sH_{g_r}'=\bigoplus_{\Sb b,c\cr b+c=|\lab|\endSb}
\bigoplus_{\Sb \be\in\La(n,b), \ga\in\La(n,c)\cr \be+\ga=\lab\endSb}
\bar x_\la T_{w_{\be,\ga}^{(2)}}e_{a+b,c}\sH_{g_r}'
\end{equation}
Moreover, we have a module isomorphism
$$x_\la\sH_{g_r}'
\cong
\bigoplus_{\Sb b,c\cr b+c=|\lab|\endSb}
\bigoplus_{\Sb \be\in\La(n,b), \ga\in\La(n,c)\cr \be+\ga=\lab\endSb}
\bar x_{(\laa,\be,\ga)}e_{a+b,c}\sH_{g_r}'.
$$
\end{cor}

\begin{pf}
The equality (\ref{lpp}) follows from Lemma \ref{lp2} and base change.
To prove the isomorphism, we use the notation introduced above to note that
$$\bar x_\la T_{w_{\be,\ga}^{(2)}}e_{a+b,c}\sH_{g_r}'=
h_{\be,\ga}T_{w_{\be,\ga}^{(2)}}\bar x_{(\laa,\be,\ga)}e_{a+b,c}\sH_{g_r}'.$$
So we have an epimorphism
from $\bar x_{(\laa,\be,\ga)}e_{a+b,c}\sH_{g_r}'$ to 
$\bar  x_{\la}T_{w_{\be,\ga}}e_{a+b,c}\sH_{g_r}'$, via
 left multiplication by $h_{\be,\ga}T_{w^{(2)}_{\be,\ga}}$.
By looking at the image of the ``standard'' basis $\{\bar x_{(\laa,\be,\ga)}e_{a+b,c}
T_w\mid w\in\bar\sD_{(\laa,\be,\ga)}\}$
, we see from Lemma \ref{lp2} that this map is
in fact an isomorphism. That is, we have
\begin{equation}\label{lppp}
\bar x_{(\laa,\be,\ga)}e_{a+b,c}\sH_{g_r}'\overset
\sim\to \bar  x_{\la}T_{w_{\be,\ga}^{(2)}}e_{a+b,c}\sH_{g_r}'.
\end{equation}
\end{pf}

\begin{thm} \label{mri}For any $\sZ_{0,g_r}$-algebra $\sZ'$, the algebras 
$\sS_q^2(r,r;\sZ')$ and
$\sS_q^{\oneb}(r,r;\sZ')$ are Morita equivalent.
In particular, this Morita equivalence holds in the linear prime case.
\end{thm}

\begin{pf} It suffices from the definition of these algebras
to prove that, for any $\la\in\Pi_2(r,r)$ with $|\laa|>0$,
every direct summand of $x_\la\sH_{g_r}'$ appearing as one of the
terms on the right in 
(\ref{lpp}) is isomorphic to a direct
summand of $x_\mu\sH_{g_r}'$ for some $\mu\in\Pi_{\oneb}(r,r)$.
To see this, we take a direct summand
$\bar  x_{\la}T_{w_{\be,\ga}^{(2)}}e_{a+b,c}\sH_{g_r}'$ of
$x_\la\sH_{g_r}'$.
By the isomorphism (\ref{lppp}), we see that the summand is isomorphic
to $\bar x_{(\laa,\be,\ga)}e_{a+b,c}\sH_{g_r}'$. 
Now, if we take $\mu=(\laa,\be,\ga)$, then, certainly,
$\bar x_\mu\sH_{g_r}'$ has a direct summand isomorphic to
$\bar x_{(\laa,\be,\ga)}e_{a+b,c}\sH_{g_r}'$.
\end{pf}

\noindent
{\bf The type $D$ case.} 
We now turn to the type $D$ case. As in \S 2, we view $\dH$ as a subalgebra
of $\sH_{1,q}$. Let $\dg_r=2\prod_{j=1}^{r-1}(q^j+1)=q^mg_r(1,q)$ and let
$\sZ_{\dg_r}$ be the ring obtained by localizing $\sZ$ at $\dg_r$.
We shall consider the algebra
$\sH_{\dg_r}=\dH\otimes_\sZ\sZ_{\dg_r}$ and define the elements 
$\dv_{a,b}=\dpi_a(T_{w_{a,b}}-T_{f(w_{a,b})})\dpi_b$,
where $f$ is the flipping map. 
Some of the properties of the elements $\dv_{a,b}$ are known
from \cite{P}. Since the proof in \cite[3.3]{P} uses dimension
comparison, we include a more general proof  below. 

\begin{lem} \label{pro} Let $a,b$ be non-negative integers.
Then:-

(a) $v_{a,b}=\cases (1+T_{s_0})\dv_{a,b},\quad\text{ if }a\ge1,\cr
\dv_{a,b}(1-T_{s_0}),\quad\text{ if }b\ge1.\cr
\endcases$

(b) $\dv_{a,b}\sH(\fS_{(a,b)})=\sH(\fS_{(b,a)}) \dv_{a,b}$.

(c)
If $a\ge1$ and $b\ge1$, then $\dv_{a,b}(T_{w_{b,a}}-T_{fw_{b,a}})\dv_{a,b}=\dv_{a,b}z$ for some 
central element $z$ of $\sH(\fS_{(b,a)})$. Moreover,  $z$
is invertible in $\dH_{\dg_r}$.

(d) $\dv_{a,r-a}\dH_{\dg_r}\cong \dv_{r-a,a}\dH_{\dg_r}$ as $\dH_{\dg_r}$-module.
\end{lem}

\begin{pf} The statement (a) follows from a direct computation
(see \cite[1.8]{HW}).
(b) is \cite[(1.14)]{P}.
To see (c), we use (a) to have 
$$\aligned
&\quad\,(1+T_{s_0})\dv_{a,b}(T_{w_{b,a}}-T_{fw_{b,a}})\dv_{a,b}\cr
&=\dv_{a,b}(1-T_{s_0})T_{b,2}(T_{s_1}-T_{s_0s_1s_0})T_{b+1,2}\cdots T_{a+b-1,a}
\dv_{a,b}\cr
&=\dv_{a,b}T_{b,2}(1-T_{s_0})T_{s_1}(1+T_{s_0})T_{b+1,2}\cdots T_{a+b-1,a}
\dv_{a,b}\cr
&=v_{a,b}T_{w_{a,b}}v_{a,b}=v_{a,b}z=
(1+T_{s_0})\dv_{a,b}z\quad \text{ by \cite[(3.23)]{DJ}},\cr
\endaligned
$$
where $T_{j,i}=T_{s_j}\cdots T_{s_i}$ for $j>i$.
Canceling $(1+T_{s_0})$, we obtain the requried relation.

We now prove (d). The result is clear if $a=0$ or $r$.
Assume now $0<a<r$ and let $b=r-a$. Then
$$\aligned
\dpi_a(T_{w_{a,b}}-T_{fw_{a,b}})\dv_{b,a}\dH_{\dg_r}
&=\dv_{a,b}(T_{w_{b,a}}-T_{fw_{b,a}})\dpi_a\dH_{\dg_r}\cr
&\supseteq \dv_{a,b}(T_{w_{b,a}}-T_{fw_{b,a}})\dv_{a,b}\dH_{\dg_r}\cr
&=\dv_{a,b}z\dH_{\dg_r}=\dv_{a,b}\dH_{\dg_r}\quad\text{ by (c).}\cr
\endaligned
$$
So the inclusion has to be an equality, and the required isomorphism follows.
\end{pf}

\begin{thm}\label{mrii} For any $\sZ_{\dg_r}$-algebra $\sZ'$, the algebras 
$\sS_q^{2.5}(r,r;\sZ')$ and
$\sS_q^{\oned}(r,r;\sZ')$  are Morita equivalent to
the algebra $\sS_q^{1.5}(n,r;\sZ')$.
In particular, this Morita equivalence holds in the linear prime case.
\end{thm}

\begin{pf} We first look at the Morita equivalence
between $\sS_q^{2.5}$ and $\sS_q^{1.5}$.
Since $\dT^{2.5}_{\sZ'}\cong \sT^2_{\sZ'}|_{\dH'}$ (\ref{clii})
and (\ref{lpp}) holds when restricted to $\dH'$, 
it suffices to
prove that any direct summand appeared in (\ref{lpp}) is isomorphic to
a direct summand of $\dT^{1.5}_{\sZ'}$. Equivalently, by (\ref{lppp}),
we only need to prove that, if $a+b=r$ and $\al$, $\be$ are compositions
of $a$, $b$ respectively, then we have an $\dH_{\dg_r}'$-module
isomorphism
\begin{equation}\label{iso}
\bar x_{(\al,\be)}e_{a,b}\dH_{\dg_r}'
\cong \bar x_{(\be,\al)}e_{b,a} \dH_{\dg_r}'.
\end{equation}
This is certainly true if $a=0$ since
$$\aligned
\bar x_{(0,\be)}e_{0,r}\dH_{\dg_r}'=
(1-T_{s_0})\bar x_{\be}\dpi_{r}\dH_{\dg_r}'&\cong
\bar x_{\be}\dpi_{0,r}\dH_{\dg_r}'\cr
&\cong
(1+T_{s_0})\bar x_{\be}\dpi_{r}\dH_{\dg_r}'=
\bar x_{(\be,0)}e_{r,0}\dH_{\dg_r}'.\cr
\endaligned$$

Assume now $a\ge1$. 
Recall from \cite[3.27]{DJ} that $e_{a,b}=v_{a,b}z_{a,b}T_{w_{b,a}}$,
where $z_{a,b}$ is central in $\sH'(\fS_{(a,b)})$, and invertible in 
$ \dH_{\dg_r}'$. Write $\ce_{a,b}=\dv_{a,b}z_{a,b}T_{w_{b,a}}$.
Then $\ce_{a,b}$ commutes with the elements of $\sH(\fS_{a,b})$
(\ref{pro}b), and
$e_{a,b}=(1+T_{s_0})\ce_{a,b}$.
Thus, 
$$\aligned
\bar x_{(\al,\be)}e_{a,b}\dH_{\dg_r}'
&=(1+T_{s_0})\bar x_{(\al,\be)}\ce_{a,b}\dH_{\dg_r}'
\cong \bar x_{(\al,\be)}\ce_{a,b}\dH_{\dg_r}'\cr
&=\bar x_{(\al,\be)}\dv_{a,b}\dH_{\dg_r}'.\cr
\endaligned
$$

We need to deal with two cases: $a=1$ and $a>1$. 
If $a=1$, then there is no $T_{s_1}$ involved in 
$\bar x_{(1,\be)}$ and so $T_{f(w_{r-1,1})}\bar x_{(1,\be)}=
\df(\bar x_{(\be,1)})T_{fw_{r-1,1}}$. Since
$\dpi_aT_{s_1}=\dpi_{a}T_{u}$ for all $a\ge 2$, we have
$$\aligned
&\quad\,\dpi_{r-1}(T_{w_{r-1,1}}-T_{fw_{r-1,1}})
\bar x_{(1,\be)}\dv_{1,r-1}\dH_{\dg_r}'\cr
&=
\bar x_{(\be,1)}\dpi_{r-1}(T_{w_{r-1,1}}-T_{fw_{r-1,1}})\dv_{1,r-1}\dH_{\dg_r}'
\cr
&=\bar x_{(\be,1)}\dv_{r-1,1}\dH_{\dg_r}',\cr
\endaligned
$$
which implies (\ref{iso}) in this case.

Finally, if $a\ge2$, then, by symmetry, we may assume that $b\ge2$.
Noting $\dpi_aT_{s_1}=\dpi_{a}T_{u}$ again,  we have
$$\aligned
&\quad\,\dpi_{b}(T_{w_{b,a}}-T_{fw_{b,a}})
\bar x_{(\al,\be)}\dv_{a,b}\dH_{\dg_r}'\cr
&=
\dpi_{b}(\bar x_{(\be,\al)}T_{w_{b,a}}-T_{fw_{b,a}}
\df(\bar x_{(\al,\be)}))\dv_{a,b}\dH_{\dg_r}'\cr 
&=
\dpi_{b}(\bar x_{(\be,\al)}T_{w_{b,a}}-
\df(\bar x_{(\be,\al)})T_{fw_{b,a}})\dv_{a,b}\dH_{\dg_r}'\cr
&=\bar x_{(\be,\al)}\dpi_{b}(T_{w_{b,a}}-T_{fw_{b,a}})\dv_{a,b}\dH_{\dg_r}'\cr
&=\bar x_{(\be,\al)}\dv_{b,a}\dH_{\dg_r}',\cr
\endaligned
$$
which implies (\ref{iso}).

We now prove
the Morita equivalence between $\sS_q^\oned$ and $\sS_q^{1.5}$.
Since $\dT^{\oned}_{\sZ'}\cong \sT^\oneb_{\sZ'}|_{\dH'}$
(\ref{onedp}(1)), the proof above together with an argument
similar to the proof of (\ref{mri}) shows that, when  
$\dT^{\oned}_{\sZ'}$ and $\dT^{1.5}_{\sZ'}$ decompose
via (\ref{lpp}), they  have the same direct summands
(up to multiplicity).
Therefore, the requried Morita equivalence follows.
\end{pf}

\begin{rems} \label{onedd}
(1) We point out that, in
the type $B$ case, results \cite[(4.17)]{DJ} and \cite[(7.6)]{GrH}
give an interpretation of the $q$-Schur$^{\oneb}$ algebra
 in terms of $q$-Schur algebras,
effectively making these algebras quasi-hereditary. However,
such a result as \cite[(4.17)]{DJ} is not available in the type $D$
case, even with the linear prime assumption, so that (\ref{mrii}) is new
in this case. Together with (\ref{straa}) below and 
 a result  \cite[7.15]{GrH} of Gruber-Hiss,
it implies that the
$q$-Schur$^{\oned}$ algebra is 
also quasi-hereditary under the linear prime assumption.
See  Remark \ref{char2}b below. 

(2) Theorem \ref{mrii} does not give a Morita equivalence
between $\sS_q^\oned$ and $\sS_q^{\onedp}$ (see (\ref{onedp}(2))).
It can be checked  that such a Morita equivalence
exists for group algebras, but with a proof involving, when $r$ is even,
 non-parabolic subgroups.
So there is no obvious ``$q$-analogue'', and we do not
know if it exists, when $r$ is even.
We point out that, without such a Morita equivalence,
result  \cite[7.17]{GrH} holds only for the $q$-Schur$^\oned$
algebra in the type $D$ case.
\end{rems}

\section{Hom spaces between $x_\la\sH'$ and $y_\mu\sH'$}

We shall assume $\sH'=\sH'_{q_0,q}$ until further notice.
(We do not assume $q_0=1$ in the next three sections.)
From Lemma \ref{6b}(a), we see that the study of the Hom space
$\Hom_{\sH'}(x_\la^\eta\sH',x_\mu \sH')$ is important. In particular
we need the base change property
$$\Hom_{\sH}(x_\la^\eta\sH,x_\mu \sH)\otimes\sZ' \cong\Hom_{\sH'}(x_\la^\eta\sH',x_\mu \sH')\cong x_\mu \sH'\cap\sH' x_\la^\eta,$$
which eventually gives the homological property Theorem \ref{hp3}
below as described in
\cite[(1.2.9-10)]{DPS1}. We need some preparation in this section and the next.
We first look at  the intersection
$\pi_a \sH'\cap\sH'\pi_b^-$.

\begin{lem}\label{piHpi1} For any non-negative integers $a,b$ with $a+b\le r$,
the set $\sB_{a,b^-}:=\{\pi_a T_x\pi_b^-\mid x\in\sD_{\om_a,\om_b}^0\}$ forms a basis
of the $\sZ'$-module $\pi_a\sH'\pi_b^-$. Moreover, we have
$\sB_{a,b^-}=\{\pi_a T_x\pi_b^-\mid x\in\sD_{\om_a,\om_b}\}\backslash\{0\}$.
(That is, every $\pi_a T_x\pi_b^-$ is either zero or a basis element.)
A similar basis $\sB_{a^-,b}$ may be constructed for 
$\pi_a^-\sH'\pi_b$.
\end{lem}

\begin{pf}
We first prove that the set generates $\pi_a\sH'\pi_b^-$. 
Clearly, $\pi_a\sH'\pi_b^-$ is generated by $\pi_a T_w\pi_b^-$ ($w\in W$).
If $w=x\hat dy$ with $x\in W_a$, $y\in W_b$ is the Howlett
decomposition,
then $\pi_a T_w\pi_b^-=\pi_a T_xT_{\hat d}T_y\pi_b^-
=\pi_a h_1T_{\hat d}h_2\pi_b^-$ for some $h_1\in\sH'(\oW_a)$
and $h_2\in\sH'(\oW_b)$. By \cite[2.2.7]{DS2}, we have
$h_1T_{\hat d}h_2=\sum_{z\in\sD_{\om_a,\om_b}}\al_zT_z$. It follows that 
$\pi_a\sH'\pi_b^-$ is generated by $\{\pi_a T_w\pi_b^-\}_{w\in\sD_{\om_a,\om_b}}$.
Suppose now that $d\in\sD_{\om_a,\om_b}\backslash\sD_{\om_a,\om_b}^0$.
 Let $d=u\hat dv$
be the right distinguished  decomposition (\cite[\S2.3]
{DS2}), where $\hat d\in\sD_{\hat\om_a,\hat\om_b}$. 
Then $s_0\in C_{\om_a}^d\cap C_{\om_b}$,
since $C_{\om_a}^d\cap C_{\om_b}=C_{\om_a}^{\hat d}\cap C_{\om_b}$.
Thus, $(1+T_{s_0})T_{\hat d}=T_{\hat d}(1+T_{s_0})$. Therefore,
$\pi_a T_d\pi_b^-=T_u\pi_a T_{\hat d}\pi_b^-T_v=0$.
This proves that $\sB_{a,b^-}$ generates $\pi_a\sH'\pi_b^-$, and
the last assertion follows.

We now prove that the elements of $\sB_{a,b^-}$ are linearly independent.
Assume $h=\sum_{x\in\sD_{\om_a,\om_b}^0}\al_d \pi_a T_d\pi_b^-=0$ for some 
$\al_d\in\sZ'$.
Since 
$$T_{t_1^{\ep_1}\cdots t_a^{\ep_a}}T_dT_{t_1^{\ep_{a+1}}\cdots 
t_b^{\ep_{a+b}}}\in
\sH'(W_a  \hat dW_b):=\sum_{w\in W_a\hat dW_b}\sZ'T_w\,\,(\ep_i\in\{0,1\}),$$
we have $0=pr_{\hat d}(h)=\sum_{x\in\sD_{\om_a,\om_b}^0\cap W_a\hat dW_b}
\al_d \pi_d T_x\pi_b^-=0$, where $pr_{\hat d}$
is the projection onto $\sH'(W_a  \hat dW_b)$. Now, by a further consideration of
the projection onto 
$\sH'(\oW_a \hat d\oW_b)$, we see that all
$\al_x=0$. This proves the linear independence, 
and hence, the lemma.
\end{pf}

\begin{cor}\label{piHpi2}  
For any non-negative integers $a,b$ with $a\ge1$, the set
$$
\aligned
\{\pi_{a-1}T_{t_a}^{\ep_1}T_{t_{a+1}}^{\ep_2}T_y\pi_b^-,
\pi_{a-1} T_{s_a}^{\ep_3}T_{t_{a}}^{\ep_4}&T_d\pi_b^-\mid\, 
y\in\sD_{\om_{a+1},\om_b}^0, \ep_j\in\{0,1\}\,\forall j,\cr
&d\in\sD_{\om_a,\om_b}^0\backslash
(\sD_{\om_{a+1},\om_b}^0\cup
t_{a+1}\sD_{\om_{a+1},\om_b}^0)\}.\endaligned
$$ is part of a basis for $\pi_{a-1}\sH'\pi_b^-$.
\end{cor}

\begin{pf} Let $t=t_a^{\ep_1}t_{a+1}^{\ep_2}$. 
We see by writing $t=gz$, where $z\in\oW_{a+1}$ and
$g\in t\oW_{a+1}$ is a distinguished left coset representative, that, for any 
$u\in \oW_{a+1}$, if $T_x$ appears as a term in $T_tT_u$ (i.e. with
non-zero coefficient), then $x=tx'$ for some $x'\in\oW_{a+1}$.
So $n_0(cx)=n_0(c)+n_0(x)$ for all $c\in C_{\om_{a-1}}$. By 
Lemma \ref{tip}c, we have
$x\in\sD_{\om_{a-1}}\cap W_{a+1}$.
For $w\in\sD_{\om_{a+1},\om_b}^0$, let 
$w=u\hat w v$ be the right distinguished decomposition with $\hat w\in
W_{a+1}wW_b$ distinguished (for parabolic subgroups). 
Recall from \cite[(2.3.1)]{DS2} that $v\in\oW_b$. Also,
since $C_{\om_{a+1}}^{\hat w}\cap C_{\om_b}=1$, the distinguished
left decomposition $u'\hat wv'$ of $x\hat w$ has
$v'\in\oW_b$. Consequently \cite[(2.2.5)]{DS2}, $x\hat w\in\sD_{\om_b}^{-1}$.
On the other hand, since
$x\in\sD_{\om_{a-1}}\cap \oW_{a+1}$, we find $x\hat w\in\sD_{\om_{a-1}}$.
(Write $x=x_1x_2$ with $x_2$ distinguished in $W_{a-1}x$, and apply
\cite[(2.2.7)]{DS2}.)
Then
$$T_tT_w=(T_tT_u)T_{\hat w}T_v=
\sum_{x\in\sD_{\om_{a-1}}\cap W_{a+1}}{\al_x}T_{x\hat w}T_v
=q^mT_{tw}+\sum_{tw<z}\be_z T_z,$$
where all $z\in\sD_{\om_{a-1},\om_b}$ by 
\cite[(2.2.7)]{DS2} again. 
Since $\pi_{a-1}T_z\pi_b^-=0$ for $z\not\in\sD_{\om_{a-1},\om_b}^0$,
we have 
$$\pi_{a-1}T_tT_w\pi_b^-=q^m\pi_{a-1}T_{tw}\pi_b^-+\sum_{tw<z,z\in\sD_{\om_{a-1},\om_b}^0}\be_z \pi_{a-1}T_z\pi_b^-.$$
A similar argument shows that this is also true for
$t=s_a^{\ep_3}t_{a}^{\ep_4}$ and $w\in
\sD_{\om_a,\om_b}^0\backslash(\sD_{\om_{a+1},\om_b}^0\cup
t_{a+1}\sD_{\om_{a+1},\om_b}^0)$.
Therefore, the given set is linearly independent and can be extended
to a basis for $\pi_{a-1}\sH'\pi_b^-$ by Cor. \ref{tip11}
and Lemma \ref{piHpi1}.
\end{pf}

We now can show that the intersection
$\pi_a \sH'\cap\sH'\pi_b^-$ is free.

\begin{thm}\label{piHpi3} Assume that $q_0+1$ is not a zero
divisor in $\sZ'$. Then,
for any positive integers $a,b\le r$, we have
$$\pi_a\sH'\cap\sH'\pi_b^-= \pi_a\sH'\pi_b^-.$$
In particular, if $a+b>r$, then the intersection is 0.
\end{thm}

\begin{pf} 
We apply induction on $a$. Note first
the elements, $T_w\pi_b^-, w\in\sD_{\om_b}^{-1}$,
form a basis of $\sH'\pi_b^-$ as in (the proof of)
\cite[(3.2.2a)]{DS2}. Now let $a=1$ and consider
$h=\sum_{w\in\sD_{\om_b}^{-1}}
\al_wT_w\pi_b^-\in \pi_1\sH'$.
We have $T_{s_0}h=q_0h$. For $w\in\sD_{\om_a}^{-1}$ write
$w=u\hat w v$ with $\hat w\in\langle s_0\rangle w W_b$ distinguished
and $u\in \langle s_0\rangle$, $v\in\oW_b$,
and put $V(w):=\langle s_0\rangle^{\hat w}\cap W_b$,
a parabolic subgroup.
Then $T_{s_0}T_{\hat wv}\pi_b^-=-T_{\hat wv}\pi_b^-$ if
$s_0\in V(w)$. Thus,
$$\aligned
T_{s_0}h
&=\sum_{s_0w>w}\al_w T_{s_0w}\pi_b^-+\sum_{s_0w<w}\al_w T_{s_0}T_w\pi_b^-\cr
&=\sum_{\Sb s_0w>w\cr s_0\in V(w)\endSb}-\al_w T_{w}\pi_b^-
+\sum_{\Sb s_0w>w\cr s_0\not\in V(w)\endSb}\al_w T_{s_0w}\pi_b^-\cr
&\quad\, +\sum_{s_0w<w}(q_0-1)\al_w T_{w}\pi_b^-+
\sum_{s_0w<w}q_0\al_w T_{s_0w}\pi_b^-\cr
\endaligned
$$
Note that the last equation gives $T_{s_0}h$ as a linear
combination of basis elements $T_y\pi_b^-, y\in\sD_{\om_b}^{-1}$.
Thus, equating  the coefficients of $T_w\pi_b^-$ in $T_{s_0}h=q_0h$,
we obtain that 
$$\cases
(1)\, \text{ if }s_0w<w, s_0\not\in V(w), &\text{ then } \al_{s_0w}+(q_0-1)\al_w=q_0\al_w;\cr
(2)\, \text{ if }s_0w<w, s_0\in V(w), &\text{ then } (q_0-1)\al_w=q_0\al_w;\cr
(3)\, \text{ if }s_0w>w ,s_0\in V(w), &\text{ then }-\al_w+q_0\al_{s_0w}=q_0\al_w.\cr
\endcases
$$
There is a fourth case, but
we do not require it. From (2), we have $\al_w=0$, for $s_0w<w, s_0\in V(w)$.
By (3) and noting that $q_0+1$ is not a zero-divisor, we have 
$\al_w=0$, for $s_0w>w, s_0\in V(w)$. Thus, 
$\al_{s_0w}=\al_w=0$ for all $w$ with $s_0\in V(w)$.
From (1) we see that $\al_{s_0w}=\al_w$ for all $w$ with  $s_0\not\in 
V(w)$.
Therefore,  $h$ can be written as $(1+T_{s_0})h'\pi_b^-=
\pi_1h'\pi_b^-$ for some $h'\in\sH'$, 
proving the case for $a=1$.

 Assume now
$a\ge 1$ and $\pi_a\sH'\cap\sH'\pi_b^-=\pi_a\sH'\pi_b^-.$
Then, $T_{s_a}\pi_a\sH'\cap\sH'\pi_b^-=T_{s_a}\pi_a\sH'\pi_b^-$.
We first prove that the result holds for $a+1$ with $a+1+b\le r$.
By \cite[(4.1.2)]{DS2}, we have $\pi_{a+1}\sH'=\pi_a\sH'\cap T_{s_a}\pi_a\sH'$.
Thus, by induction, we obtain
\begin{equation}\label{php}
\pi_{a+1}\sH'\cap\sH'\pi_b^-=\pi_a\sH'\cap T_{s_a}\pi_a\sH'\cap\sH'\pi_b^-
=\pi_a\sH'\pi_b^-\cap T_{s_a}\pi_a\sH'\pi_b^-.
\end{equation}
We need to prove that
$$\pi_a\sH'\pi_b^-\cap T_{s_a}\pi_a\sH'\pi_b^-=\pi_{a+1}\sH'\pi_b^-.$$
The inclusion ``$\supseteq$'' is clear.
Suppose $h\in \pi_a\sH'\pi_b^-\cap T_{s_a}\pi_a\sH'\pi_b^-$. Then,
by Lemma \ref{piHpi1},
$h=\pi_ah_1\pi_b^-= T_{s_a}\pi_ah_2\pi_b^-$
where $h_1=\sum_{x\in\sD_{\om_a,\om_b}^0}\al_x T_x$ and
$h_2=\sum_{y\in\sD_{\om_a,\om_b}^0}\be_y T_y$. So we have
$(\pi_ah_1- T_{s_a}\pi_ah_2)\pi_b^-=0$. Since 
$\pi_{a+1}\sH'=\pi_a\sH'\cap T_{s_a}\pi_a\sH'$, if we could prove that
$h_0=\pi_ah_1- T_{s_a}\pi_ah_2=0$, i.e., $h'=\pi_ah_1=T_{s_a}\pi_ah_2$,
then we would have $h'\in \pi_{a+1}\sH'$, and therefore, 
$h=h'\pi_b^-\in\pi_{a+1}\sH'\pi_b^-$,
as desired. So it remains to prove $h_0=0$.
For simplicity of notation,
we put $D=\sD_{\om_a,\om_b}^0\backslash( \sD_{\om_{a+1},\om_b}^0
\cup
t_{a+1}\sD_{\om_{a+1},\om_b}^0)$
and $B=B_1\cup B_2$ where
 $$\aligned
B_1&=\{T_x, T_{t_{a}}T_x, T_{t_{a+1}}T_x,T_{t_{a}}T_{t_{a+1}}T_x\mid
x\in\sD_{\om_{a+1},\om_b}^0\},\text{ and }\cr
B_2&=\{T_d, T_{t_{a}}T_d, T_{s_{a}}T_d,T_{s_{a}}T_{t_{a}}T_d\mid
d\in D\}.
\endaligned
$$
Since $a>0$, Cor. \ref{piHpi2} gives 
the linear independent set  $\sB=\pi_{a-1}B\pi_b^-$
in $\pi_{a-1}\sH'\pi_b^-$. 
Write
$$T_{s_a}\pi_ah_2=\pi_{a-1}T_{s_a}(q^{a-1}+T_{t_a})h_2
=\pi_{a-1}(q^{a-1}T_{s_a}h_2+T_{s_a}T_{t_a}h_2),$$
and
$$h_2=\sum_{y\in\sD_{\om_{a+1},\om_b}^0}\be_y T_y+
\sum_{y\in\sD_{\om_{a+1},\om_b}^0}\be_{t_{a+1}y} T_{t_{a+1}}T_y
+\sum_{d\in D}\be_dT_d.$$
Since $y\in\sD_{\om_{a+1},\om_b}^0$ implies
 $s_ay,s_at_{a+1}y, t_{a+1}y\in
\{1,t_a,t_{a+1}\}\sD_{\om_{a+1},\om_b}^0$
(see Lemma \ref{tip}b), and $d\in D$ implies 
$s_ad\in s_aD$, 
we see that $T_{s_a}h_2$ is a linear
combination of the elements in $B$. Likewise,
$T_{s_a}T_{t_{a}}T_y$ (resp. $T_{s_a}T_{t_{a}}T_{t_{a+1}}T_y$)
is a linear combination of
 $T_{t_{a+1}}T_w$ (resp. $T_{t_{a}}T_{t_{a+1}}T_w$)
with $w\in\sD_{\om_{a+1},\om_b}^0$, and $T_{s_a}T_{t_a}T_d$ is
a linear combination of the elements in $B_2$.
Hence $T_{s_a}T_{t_a}h_2$
is a linear combination of the elements of $B$.
On the other hand,
$\pi_ah_1=\pi_{a-1}(q^{a-1}+T_{t_a})h_1$ and
$(q^{a-1}+T_{t_a})h_1$ is a linear
combination of the elements of $B$.
Therefore, we see that $h_0$ is a linear
combination of $\pi_{a-1}T_w$ with $w\in B$. 
Since  $\sB=\pi_{a-1}B\pi_b^-$ is a linear independent set,
it follows from $h_0\pi_b^-=0$ that $h_0=0$. 

Assume now that $a+b\ge r$. Then $\pi_{a+1}\sH'\pi_b^-=0$ by \cite[3.2]{DJ}.
We need to prove that $\pi_{a+1}\sH'\cap\sH'\pi_b^-=0$.
By (\ref{php}), it suffices to prove that $\pi_a\sH'\pi_b^-\cap T_{s_a}\pi_a\sH'\pi_b^-=0$ if $a+b=r$.
Let $w_{a,b}$ be the (unique) distinguished $\oW_a$-$\oW_b$ double
coset representative with the trivial intersection property:
 $\oW_a^{w_{a,b}}\cap \oW_b=\{1\}$, 
as given at the beginning of \S 3,
and let 
$\pi(a,b^-)=\pi_a T_{w_{a,b}}\pi_b^-$
(cf. footnote 5 below). Then, by \cite[3.11]{DJ},
$\pi_a\sH'\pi_b^-=\sH'(\oW_{(a,b)})\pi(a,b^-)$.
Now the linear independence (see  \cite[3.15]{DJ}) between bases
$\sB=\{T_w\pi(a,b^-)\mid w\in \oW_{(a,b)}\}$ for $\pi_a\sH'\pi_b^-$
and $T_{s_a}\sB$ for $T_{s_a}\pi_a\sH'\pi_b^-$ implies that 
the intersection is trivial. This completes the proof of the theorem.
\end{pf}

We have the following generalization which has been known for 
parabolic subgroups.
Recall the definition of $y_\la$ in Remark \ref{djm}c.

\begin{thm}\label{piHpi6} Assume that both $q_0+1$ and $q+1$ are 
not zero divisors in $\sZ'$. For bicompositions $\la,\mu$ of $r$,
we have $x_\la\sH'\cap\sH'y_\mu=x_\la\sH'y_\mu.$ Moreover,
this intersection is $\sZ'$-free with basis $\{x_\la T_w y_\mu\mid w\in\sD_{\la,\mu}^0\}$.
A similar result holds for $y_\la\sH'\cap\sH'x_\mu$.
\end{thm}

\begin{pf}  Clearly, we have $x_\la\sH'\cap\sH'y_\mu
\supseteq x_\la\sH'y_\mu.$
Applying \cite[(4.1.3)]{DS2}\footnote{If we put $q_{s_0}=q_0$, 
$q_{s_i}=q$ for $1\le i\le r-1$,
and define for $w\in W$, $q_w=q_{s_{i_1}}\cdots q_{s_{i_m}}$,
 where
$w=s_{i_1}\cdots s_{i_m}$ is a  reduced expression,
then the {\it ring} homomorphism on $\sH'$ sending
$q_s$ to $q_s^{-1}$ and $T_w$ to $(-1)^{\ell(w)}q_wT_w$ will
interchange $x_\la$ with $y_\la$. Thus, \cite[(4.1.3)]{DS2}
holds for $y_\la\sH'$, and also for $\sH'y_\la$.}
 and Theorem \ref{piHpi3}, we obtain
$x_\la\sH'\cap\sH'y_\mu=\bar x_\la\sH'\cap \sH'\bar y_\mu\cap \pi_a\sH'\pi_b^-$. By \cite[1.1h]{DPS3}, we have $\bar x_\la\sH'\cap \sH'\bar y_\mu=
\bar x_\la\sH'\bar y_\mu$ since both $\oW_{\la}$ and $\oW_{\mu}$ are
parabolic. (This requires that $q+1$ is not a zero divisor.)
Suppose $h\in\bar x_\la\sH'\bar y_\mu\cap \pi_a\sH'\pi_b^-$ and write
$h=\sum_{w\in\sD_{\om_a,\om_b}^0}\al_w\pi_a T_w\pi_b^-$. Then,
$T_sh=qh$ and
$hT_t=-h$ for all $s=s_i\in\oW_\la$ and $t=s_j\in\oW_\mu$.
Since we have  $sw, wt\in\sD_{\om_a,\om_b}^0$
whenever $w\in\sD_{\om_a,\om_b}^0$ (see Lemma \ref{tip}b), it follows
from Lemma \ref{piHpi1},
for $h'=\sum_{w\in\sD_{\om_a,\om_b}^0}\al_w T_w$ with $h=\pi_a h'\pi_b^-$,
that $T_sh'=qh'$ and
$h'T_t=-h'$ for all $s=s_i\in\oW_\la$ and $t=s_j\in\oW_\mu$.
Therefore, $h'=\sum_{w\in\sD_{\om_a,\om_b}^0\cap\sD_{\bar\la,\bar\mu}^0}\al_w \bar x_\la T_w
\bar y_\mu$ and 
$h=\sum_{w\in\sD_{\bar\la,\bar\mu}^0\cap\sD_{\om_a,\om_b}^0}\al_w x_\la T_w y_\mu,$ proving the first assertion.
 Since $\sD_{\bar\la,\bar\mu}^0\cap\sD_{\om_a,\om_b}^0=\sD_{\la,\mu}^0$,
the basis assertion follows immediately.
The final claim may be established by applying the standard anti-automorphism
$T_w\mapsto T_{w^{-1}}$.
\end{pf}
 
\begin{cor}\label{piHpi5} Maintain the notation above and 
 assume that both $q_0+1$ and $q+1$ are 
not zero divisors in $\sZ'$. Then
$\Hom_{\sH'}(y_\mu\sH',x_\la\sH')$ is free and
$$\Hom_{\sH'}(y_\mu\sH',x_\la\sH')\cong x_\la\sH'y_\mu.$$ 
A similar statement holds with the roles of $x_\la$ and $y_\mu$
interchanged.
\end{cor}

\section{Induced bistandard bases}

{\it Throughout the section, we assume $a+b=r$ and $q_0+1$ is not a zero-divisor
in $\sZ'$.}
Recall from the proof of Theorem \ref{piHpi3} that  $w_{a,b}$ is the (unique) distinguished $\oW_a$-$\oW_b$ double
coset representative with the trivial intersection property.

\begin{lem}\label{tip4} If $a+b=r$, then
$\sD_{\om_a,\om_b}^0=\oW_a {w_{a,b}}\oW_b$.
\end{lem}

\begin{pf}
Let $w\in\sD_{\om_a,\om_b}^0$, and let $w=udv$ be the right distinguished decomposition
of $w$ with $d\in W_a wW_b$ distinguished. Since $C_{\om_a}^d\cap C_{\om_b}=\{1\}$,
it follows that $n_0(t_1\cdots t_adt_1\cdots
t_b)=a+n_0(d)+b\ge r$, and hence $n_0(d)=0$ and $d\in\oW$. 
Since $d$ is distinguished and  
$$\{(1)d,\cdots (a)d\}=\{b+1,\cdots r\},$$ 
we must have $d=w_{a,b}$.
The rest of the proof follows immediately from
 the trivial intersection property of $w_{a,b}$.
\end{pf}

Put $\pi(a,b^-)=\pi_aT_{w_{a,b}}\pi_b^-$  and 
$\pi(a^-,b)=\pi_a^-T_{w_{a,b}}\pi_b$,\footnote{We use superscript $^-$
to indicate the ``minus'' part. The notation $\pi(a,b^-)$ is denoted
by $v_{a,b}$ in \cite{DJ}.}
 and write $\bar\sH'_m=\sH'(\oW_m)$.
The latter may be regarded as the image of $\sH'$ under the
obvious ($\sZ'$-module) epimorphism $h\mapsto \bar h$ from
$\sH'$ to $\bar\sH'$.

\begin{cor}\label{bis} If $a+b=r$, then $\pi_a\sH'\pi_b^-$ has basis
$\{T_u\pi(a,b^-)T_v\mid u\in\oW_a,v\in\oW_b\}$.
A similar result holds for $\pi_a^-\sH'\pi_b$. In particular, we have
$\pi_a\sH'\pi_b^-=\bar\sH'_a\pi(a,b^-)\bar\sH'_b$ and
 $\pi_a^-\sH'\pi_b=\bar\sH'_a\pi(a^-,b)\bar\sH'_b$.
\end{cor}
\begin{pf} Immediately follows from Lemma \ref{piHpi1} and the lemma above.
\end{pf}

For partitions $\al,\be$ of $a$, let $\bar\sH^{\prime\al\be}_a=\bar x_\al\bar\sH_a'
\cap\bar\sH'_a\bar x_\be$. This is a free $\sZ'$-module.

\begin{thm}  Let $\la,\mu$ be bicompositions  of $r$ such that
$|\laa|+|\mu^{(1)}|=r$. Then
$$x_\la^\eta\sH'\cap\sH'x_\mu=\bar\sH^{\prime\laa\mu^{(2)}}_a\pi(a^-,b)
\bar\sH^{\prime\lab\mu^{(1)}}_b,$$
where $a=|\laa|$ and $b=|\mu^{(1)}|$.
\end{thm}

\begin{pf} We first note that $w_{a,b}\in \bar W_{(a,b)}w_{a,b}\bar W_{(b,a)}$ is distinguished and
$\oW_{(a,b)}^{w_{a,b}}=\oW_{(b,a)}$. Thus, we have
$$\bar x_{\la}\sH'(\oW_{(a,b)})\pi(a^-,b)=
\bar x_\laa\bar \sH'_a\pi(a^-,b)\bar x_\lab\bar\sH'_b=
\pi(a^-,b)\bar x_{\la^\star}\sH'(\oW_{(b,a)})
$$
for all $\la$ with $|\laa|=a$ and $|\lab|=b=r-a$.
Thus, by \cite[(4.1.3)]{DS2} and Thereom \ref{piHpi6}
(applying $\eta$ if necessary), we have 
\begin{equation}\label{xhpi}
\aligned
x_\la^\eta\sH'\cap\sH'x_\mu
&=(x_\la^\eta\sH'\cap\pi_a^-\sH')\cap(\sH'x_\mu\cap \sH'\pi_b)\cr
&=(x_\la^\eta\sH'\cap \sH'\pi_b)\cap(\pi_a^-\sH'\cap\sH'x_\mu)
= x_\la^\eta\sH'\pi_b\cap \pi_a^-\sH'x_\mu\cr
&=(\bar x_\laa\bar \sH_a'\pi(a^-,b)\bar x_\lab\bar \sH'_b)\cap
  (\bar\sH'_a\bar x_{\mu^{(2)}}\pi(a^-,b)\bar\sH'_b \bar x_{\mu^{(1)}})\cr
\endaligned
\end{equation}
By Corollary \ref{bis}, we have clearly that, if $h_1\pi(a^-,b) h_2=
h_1'\pi(a^-,b) h_2'$, where $h_1\in \bar x_\laa\bar \sH_a'$,
$h_1'\in \bar\sH'_a\bar x_{\mu^{(2)}}$,
$h_2\in\bar x_\lab\bar \sH'_b$  and $h_2'\in \bar\sH'_b \bar x_{\mu^{(1)}}$,
then $h_1=h_1'$ and $h_2=h_2'$. Thus, $h_1\in \bar\sH^{\prime\laa\mu^{(2)}}_a$,
$h_2\in \bar\sH^{\prime\lab\mu^{(1)}}_b$, and
consequently,
we have $x_\la^\eta\sH'\cap\sH'x_\mu=\bar\sH^{\prime\laa\mu^{(2)}}_a\pi(a^-,b)
\bar\sH^{\prime\lab\mu^{(1)}}_b,$ proving the theorem.
\end{pf}

The theorem above guarantees the existence of a Murphy-type basis 
for the intersection $x_\la^\eta\sH'\cap\sH'x_\mu$. 
We are now going to describe this basis.

Recall from \cite{M} 
 the Murphy basis  $\{x_{\bfs\bft}^a:=T_{\de(\bfs)^{-1}} \bar x_\be 
T_{\de(\bft)}\}$ of $\bar\sH'_a$, where
$\bfs,\bft$ are standard $\be$-tableaux for all partitions $\be$ of $a$
and $\de(\bfu)$ is a distinguished right $\fS_\be$-coset representative
defined by $\bfu$. Let $\bar\bT^s(\be)$ denote the
set of standard $\be$-tableaux, and let,
for $\al\in\La(n,r)$, $\bar\fT^{ss}(\be,\al)$ be the set of semi-standard
$\be$-tableaux of type $\al$. These have all been generalized
to bipartitions and bicompositions in \cite[\S1]{DS2}.
Thus, for $\la\in\Pi_2$ and $\mu\in\Pi_2^+$,
the corresponding sets are denoted by $\bT^s(\mu)$ and $\fT^{ss}(\mu,\la)$,
respectively. We refer the reader to \cite[\S1]{DS2} for the notion of 
standard and semi-standard bitableaux. 
Recall also from \cite[\S1]{DS2} the $\de$-function which 
takes a standard (bi)tableau or a 
semi-standard (bi)tableau to a certain distinguished coset representative
in $\oW$ and the function $\bar{\frak f}$ (resp. $\frak f$)
 from $\bar\bT^s(\be)$ (resp. $\bT^s(\mu)$)
to $\bar\fT^{ss}(\be,\al)$ (resp. $\fT^{ss}(\mu,\la)$). We shall denote
by $\bar\bT_\fks$ (resp. $\bT_\fks$) the inverse image of $\fks$ under 
$\bar{\frak f}$ (resp.
$\frak f$).

Murphy basis induces a basis 
(see \cite{M} and \cite{DS2}) $m^a_{\fks\fkt}$ for
 the module $\bar\sH_a^{\prime\al\be}$, for any given compositions 
$\al,\be$ of $a$,
where $\fks\in\bar\fT^{ss}(\ga,\al)$, $\fkt\in\bar\fT^{ss}(\ga,\be)$, and
$\ga$ is a partition of $a$.
Note that the basis element is simply a sum of certain Murphy 
basis elements:
$$m^a_{\fks\fkt}=\sum_{\bfs\in\bar\bT_\fks,\bft\in\bar\bT_\fkt}x_{\bfs\bft}^a.$$
We have immediately the following.

\begin{thm}\label{murphy} Let $\la,\mu$ be bicompositions  of $r$ such that
$|\laa|+|\mu^{(1)}|=r$. Then the set
$\sY_{\la^-,\mu}=\{m_{\fks_1\fkt_2}^a\pi(a^-,b) m_{\fks_2\fkt_1}^b\}$,
where $a=|\laa|$, $b=|\mu^{(1)}|$ and
 $\{m_{\fks_1\fkt_2}^a\}$ and $\{m_{\fks_2\fkt_1}^b\}$ are  bases for 
$\bar\sH^{\prime\laa\mu^{(2)}}_a$ and $\bar\sH^{\prime\lab\mu^{(1)}}_b$, 
respectively,
forms a basis for $x_\la^\eta\sH'\cap\sH'x_\mu$.
In particular, the Murphy bases for $\bar\sH'_a$ and $\bar\sH'_b$
induce a Murphy type basis
$\sY_{a^-,b}=\sY_{\om_a^-,\om_b}$
 for $\pi_a^-\sH'\pi_b=\bar\sH'_a\pi(a^-,b)\bar\sH'_b$. 
\end{thm}

The bases $\sY_{\la^-,\mu}$ for $x_\la^\eta\sH'\cap\sH' x_\mu$
 are constructed by using the Murphy type bases
for the type $A$ intersections $\bar\sH^{\prime\al\be}$. However,
we are also able to construct them directly by using 
Murphy type bases $\sX_{\la,\mu}$, defined in 
\cite{DS2} and called the bistandard bases, for the  type $B$ intersections
$x_\la\sH'\cap\sH' x_\mu$. We will see that such a 
construction  gives us certain homological properties required
in stratifying an endomorphism algebra.

Let $\sY_{a,b^-}=\eta(\sY_{a^-,b})$ and
$\Pi_{2,a}=\{\la\in\Pi_2\mid |\laa|=a\}$. Then
$$\sY_{a,b^-}=\{x_{\bfs_1\bft_1}^a\pi(a,b^-) x_{\bfs_2\bft_2}^b\mid
\bfs_i,\bft_i\in\bar \bT^s(\nu^{(i)}),
\nu\in\Pi_{2,a}^+\}.$$
replacing $\pi(a,b^-)$ by $\pi(a^-,b)$, we obtain
$\sY_{a^-,b}$. 
We shall call the bases $\sY_{a,b^-}$ and $\sY_{a^-,b}$
{\it induced bistandard} bases.

Recall from \cite[\S5.2]{DS2}  the bistandard basis of $
x_\la\sH'\cap\sH' x_\mu$:
$$\sX_{\la,\mu}=\{X_{\fks\fkt}^{\nu}\mid \nu\in\Pi^+_2,\fks\in\fT^{ss}(\nu,\la),
\fkt\in\fT^{ss}(\nu,\mu)\},$$
where $X_{\fks\fkt}^\nu=\sum_{\bfs\in\bT_\fks,\bft\in\bT_\fkt}T_{\de(\bfs)^{-1}}x_\nu T_{\de(\bft)}.$
Since $\fT^{ss}(\nu,\om_0)=\bT^s(\nu)$ where $\om_0=(-,1^r)$, the set
$$\sX_{\la,\om_0}=\{X_{\fks\bft}^{\nu}\mid \nu\in\Pi^+_2,
\fks\in\fT^{ss}(\nu,\la),
\bft\in\bT^s(\nu)\},$$
respectively,
$$\sX_{\om_0,\mu}=\{X_{\bfs\fkt}^{\nu}\mid \nu\in\Pi^+_2,
\fkt\in\fT^{ss}(\nu,\mu),
\bfs\in\bT^s(\nu)\}$$
is a basis for $x_\la\sH'$, respectively, for $\sH'x_\mu$.

\begin{thm} \label{ind} Assume $a+b=r$.
We have that
$\sY_{a^-,b}=\pi_a^-\sX_{\om_0,\om_b}\backslash\{0\}$
and
$\sY_{a,b^-}=\sX_{\om_a,\om_0}\pi_b^-\backslash\{0\}$
form bases for  $\pi_a^-\sH'\pi_b$ and $\pi_a\sH'\pi_b^-$,
respectively.
\end{thm}

\begin{pf}  We prove the first case. The second case 
can be obtained by applying the involution
$\iota:T_w\mapsto T_{w^{-1}}$ to the first. Since 
$X_{\bfs\fkt}^\nu=\sum_{\bft\in\bT_\fkt}T_{\de(\bfs)^{-1}}x_\nu T_{\de(\bft)},$
 we have by Lemma \ref{piHpi1} 
$\pi_a^-X_{\bfs\fkt}^\nu=0$ whenever 
$|\nu^{(1)}|>b=r-a$, or $|\nu^{(1)}|=b$ and $\de(\bfs)^{-1}\not\in\sD_{\om_a,\om_b}^0$.
Assume now $\pi_a^-X_{\bfs\fkt}^\nu\neq0$. Then 
$|\nu^{(1)}|=b$ and $\de(\bfs)^{-1}\in\sD_{\om_a,\om_b}^0$.

We now note that, if $\fkt=(\fkt_1,\fkt_2)\in \fT^{ss}(\nu,\om_b)$,
then, by the definition \cite[(1.2.2)]{DS2}, $\fkt_1$ is semi-standard and
contains a $\nu^{(1)}$-tableau of type $(1^b)$. Thus, $|\nu^{(1)}|=b$ forces
that $\fkt_1$ is a standard $\nu^{(1)}$-tableau, i.e.,  $\fkt_1\in
\bar\bT^s(\nu^{(1)})$. Therefore, 
 the set $\fT^{ss}(\nu,\om_b)$ can be identified
with $\bar\bT^s(\nu^{(1)})\times\bar\bT^{s}(\nu^{(2)})w_{a,b}$.
So $|\bT_\fkt|=1$ and, for $\fkt\in\fT^{ss}(\nu,\om_b),
\bfs\in\bT^s(\nu)$,
$X_{\bfs\fkt}^\nu=\sum_{\bft\in T_\fkt}X_{\bfs\bft}^\nu=X_{\bfs\bft}^\nu$
where $\fkt=\{\bft\}=\{(\bft_1,\bft_2 w_{a,b})\}$. On the other hand,
since $\sD_{\om_a,\om_b}^0=\oW_a w_{a,b}\oW_b$ by Lemma \ref{tip4},
it follows that  $\de(\bfs)=\de(\bfs_1)w_{b,a}\de(\bfs_2)$ for some 
$(\bfs_1,\bfs_2)\in\bar\bT^s(\nu^{(1)})\times\bar\bT^s(\nu^{(2)})$.
Therefore, we have
$$\aligned
\pi_a^-X_{\bfs\fkt}^\nu
&=\pi_a^-T_{\de(\bfs)^{-1}}\pi_b\bar x_\nu T_{\de(\bft)}
=\pi_a^-T_{\bar\de(\bfs_2)^{-1}}T_{w_{a,b}}T_{\bar\de(\bfs_1)^{-1}}
\pi_b\bar x_\nu T_{\de(\bft)}\cr
&=(T_{\de(\bfs_2)^{-1}}\bar x_{\nu^{(2)}}T_{\de(\bft_2)})\pi(a,b^-)
(T_{\de(\bfs_1)^{-1}}\bar x_{\nu^{(1)}}T_{\de(\bft_1)})\cr
&=x^a_{\bfs_2\bft_2}\pi(a,b^-)x_{\bfs_1\bft_1}^b,\cr
\endaligned$$
proving the inclusion
 $\pi_a^-\sX_{\om_0,\om_b}\backslash\{0\}\subseteq\sY_{a^-,b}$.
The above relation also proves  that $\sY_{a^-,b}$ is a subset
of $\pi_a^-\sX_{\om_0,\om_b}$, hence they are equal.
\end{pf}

For $\nu\in\Pi_{2,a}^+$, let $\bT^s_{a}(\nu)=
\bar\bT^s(\nu^{(1)})\times\bar\bT^s(\nu^{(2)})w_{b,a}$. Then $\bT^s_{a}(\nu)$ is the subset
of $\bT^s(\nu)$ consisting of all standard $\nu$-bitableau 
$\bft=(\bft_1,\bft_2)$
such that $\bft_1$ has entries $1,2,\cdots,a$.

We put, for $\la \in\Pi_{2,a}$,  
$$\fT^{ss}_a(\nu,\la)=\{\fks\in\fT^{ss}(\nu,\la)\mid\bT_\fks\subseteq
\bT^s_{a}(\nu)\}={\frak f}(\bT^s_{a}(\nu)).$$

\begin{cor} \label{ind1}
We have for $\la\in\Pi_{2,a}$ and $\mu\in\Pi_{2,b}$
$$\aligned
\sY_{\om_a^-,\mu}&=\pi_a^-\sX_{\om_0,\mu}\backslash\{0\}\cr
&=\{\pi_a^-X_{\bfs\fkt}^\nu\mid \nu\in\Pi^+_{2,b},\fkt\in\fT^{ss}_b(\nu,\mu),\bfs\in\bT^s(\nu),
\de(\bfs)\in\sD_{\om_b,\om_a}^0\},\text{ and }\cr
\sY_{\la^-,\om_b}&=\eta(\sX_{\la,\om_b})\pi_b\backslash\{0\}\cr
&=
\{\eta(X_{\fks\bft}^\nu)\pi_b\mid \nu\in\Pi^+_{2,a},\fks\in\fT^{ss}_a(\nu,\la),\bft\in\bT^s(\nu),
\de(\bft)\in\sD_{\om_a,\om_b}^0\}.\cr
\endaligned$$
\end{cor}

\begin{pf} First, using an argument similar to (\ref{xhpi}), we see that $
\sY_{\om_a^-,\mu}\subseteq\pi_a^-\sX_{\om_0,\mu}$.
Suppose $X_{\bfs\fkt}^\nu\in\sX_{\om_0,\mu}$. Then,
by a similar argument
as above,
$\pi_a^-X_{\bfs\fkt}^\nu\neq0$ implies that $\de(\bfs)^{-1}\in\sD^0_{\om_a,\om_b}$
and $\nu\in\Pi_{2,b}^+$. Thus, $\de(\bfs)^{-1}=\de(\bfs_2)^{-1}w_{a,b}
\de(\bfs_1)^{-1}$ for some $\bfs_i\in\bar\bT(\nu^{(i)})$.
On the other hand, since $\nu\in\Pi_{2,b}^+$, we have $\bT_\fkt
\subseteq\bT_b^s(\nu)$ and so $\fkt\in\fT^{ss}_b(\nu,\mu)$. Therefore, by (\ref{xhpi}),
one sees easily that $\pi_a^-X_{\bfs\fkt}^\nu\in \sY_{\om_a^-,\mu}$.
\end{pf}

For {\it any} bitableau $\bfu=(\bfu_1,\bfu_2)$, 
define $\bfu^\star=(\bfu_2,\bfu_1)$. Then, ${\,}^\star$ sends a standard $\nu$-tableau
to a standard $\nu^\star$-tableau.

\begin{lem}
If $\nu\in\Pi_{2,b}^+$ and $\bft\in\bT_b^s(\nu)$, then
$\de(\bft^\star)=w_{a,b}\de(\bft)\in \sD_{\om_a,\om_b}^0$, and
$$\{w_{b,a}\de(\bft)\mid \bft\in\bT_b^s(\nu^\star)\}=
\{\de(\bfs)\mid \bfs\in 
\bT^s(\nu), \de(\bfs)\in\sD_{\om_b,\om_a}^0\}
.$$
In particular, we have
$$\bT_b^s(\nu^\star)^\star=\{\bfs\mid \bfs\in 
\bT^s(\nu), \de(\bfs)\in\sD_{\om_b,\om_a}^0\}.$$
\end{lem}

\begin{pf} Write $\bft=(\bft_1,\bft_2w_{a,b})$ for some
$\bft_i\in\bar\bT^s(\nu^{(i)})$. Then 
$\de(\bft)=\bar\de(\bft_1)(w_{b,a}\bar\de(\bft_2)w_{a,b})$ and
$\bft^\star=(\bft_2w_{a,b},\bft_1)$.
If $\bft^\nu$ denotes the standard $\nu$-bitableau
in which the numbers $1,2,\cdots, r$ appear in the same order
down successive rows in the first diagram of $\la$
and then in the second diagram, then 
$$\aligned
\bft^{\nu^\star}\de(\bft^\star)&=
(\bft_2w_{a,b},\bft_1)=(\bft_2,\bft_1w_{b,a})w_{a,b}\cr
&=(\bft^{\nu^\star}_1,\bft^{\nu^\star}_2)(\bar\de(\bft_2)
w_{a,b}\bar\de(\bft_1)w_{b,a})w_{a,b}
=\bft^{\nu^\star}\bar\de(\bft_2)w_{a,b}\bar\de(\bft_1).\cr
\endaligned$$
So $\de(\bft^\star)=w_{a,b}\de(\bft)$. 

Now the inclusion ``$\subseteq$'' follows immediately with $\nu$ replaced
by $\nu^\star$.
Conversely, suppose $\bfs\in \bT^s(\nu)$ and $\de(\bfs)\in\sD_{\om_b,\om_a}^0$. Then $\de(\bfs)=\bar\de(\bfs_1)w_{b,a}\bar\de(\bfs_2)$ for some 
$\bfs_i\in\bar\bT(\nu^{(i)})$, and 
$\bfs=\bft^\nu \de(\bfs)=(\bfs_1w_{b,a},\bfs_2)$.
So $\bfs^\star=(\bfs_2,\bfs_1 w_{b,a})\in\bT_b^s(\nu^\star)$ and
$\de(\bfs)=w_{b,a}\de(\bfs^\star)$, proving the equality.
\end{pf}

We now have an alternate description of the basis $\sY_{\la^-,\mu}$
defined in \ref{murphy}.
For $\nu\in\Pi^+_{2,b}$, $\fks^\star\in\fT^{ss}_b(\nu,\la^\star)$ and $\fkt\in\fT^{ss}_b(\nu,\mu)$,
let 
$$Z_{\fks\fkt}^\nu=\sum_{\bfs\in\bT_\fks}\pi_a^-X_{\bfs\fkt}^\nu.$$

\begin{thm} \label{ind2}
We have, for $\la\in\Pi_{2,a},\mu\in\Pi_{2,b}$, the set
$$\{Z_{\fks\fkt}^\nu\mid \nu\in\Pi^+_{2,b},
\fks^\star\in\fT^{ss}_b(\nu,\la^\star),\fkt\in\fT^{ss}_b(\nu,\mu)\}$$
forms a basis for $x_\la^\eta\sH'\cap\sH'x_\mu$.
Moreover, it coincides with the basis
$\sY_{\la^-,\mu}$ defined in (\ref{murphy}).
\end{thm}

\begin{pf} Note first that $\fks\in\fT_b(\nu^\star,\la)$ is
a $\nu^\star$-tableau of type $\la$ (not necessarily semi-standard)
and that,
if
$\bft\in\bT_{\fkt}$ with $\fkt\in\fT^{ss}_b(\nu,\mu)$, then 
$\bft=(\bft_1,\bft_2w_{a,b})$ and
$\de(\bft)=\bar\de(\bft_1)(w_{b,a}\bar\de(\bft_2)w_{a,b})$ 
for some 
$\bft_i\in\bar\bT(\nu^{(i)})$,  where
$\bar\de(\bft_1)
\in\fS_{\{1,\cdots,b\}}$ and $
w_{a,b}^{-1}\bar\de(\bft_2)w_{a,b}
\in\fS_{\{b+1,\cdots,r\}}$.
Thus, $\bft^\star=(\bft_2w_{a,b},\bft_1)\in\bT_{\fkt^\star}$
and $\fkt^\star\in\fT_a(\nu^\star,\mu^\star)$.
So we have
$$\aligned
Z_{\fks\fkt}^\nu&=
\sum_{\bfs\in\bT_\fks,\bft\in\bT_{\fkt}}\pi_a^-T_{\de(\bfs)^{-1}}\pi_b\bar 
x_\nu T_{\de(\bft)}
\cr
&=\sum_{\bfs\in\bT_\fks,\bft\in\bT_{\fkt}}\pi_a^-
T_{\bar\de(\bfs_2)^{-1}}T_{w_{a,b}\bar\de(\bfs_1)^{-1}}\pi_b\bar 
x_\nu T_{w_{b,a}\bar\de(\bft_2)w_{a,b}}T_{\bar\de(\bft_1)}\cr 
&=\sum_{\bfs^\star\in\bT_{\fks^\star},\bft\in\bT_{\fkt}} 
T_{\bar\de(\bfs_2)^{-1}}T_{w_{a,b}\bar\de(\bfs_1)^{-1}w_{b,a}}\bar x_{\nu^\star}
(\pi_a^-T_{w_{a,b}}\pi_b)T_{w_{b,a}\bar\de(\bft_2)w_{a,b}}T_{\bar\de(\bft_1)}\cr
&=\sum_{\bfs^\star\in\bT_{\fks^\star},\bft^\star\in\bT_{\fkt^\star}} 
T_{\bar\de(\bfs_2)^{-1}}T_{w_{a,b}\bar\de(\bfs_1)^{-1}w_{b,a}}\bar x_{\nu^\star}
\pi_a^-T_{\bar\de(\bft_2)w_{a,b}}T_{\bar\de(\bft_1)}\pi_b\cr
&=\sum_{\bfs^\star\in\bT_{\fks^\star},\bft^\star\in\bT_{\fkt^\star}}
T_{\de(\bfs^\star)^{-1}}\bar
x_{\nu^\star}\pi_a^-T_{\de(\bft^\star)}\pi_b\cr
&=\sum_{\bft^\star\in\bT_{\fkt^\star}}\eta(X^{\nu^\star}_{\fks^\star,\bft^\star})\pi_b\cr
\endaligned$$
So $Z_{\fks\fkt}^\nu\in x_\la^\eta\sH'\cap\sH'x_\mu$ since
$x_\la^\eta\sH'\cap\sH'x_\mu= x_\la^\eta\sH'\pi_b\cap \pi_a^-\sH'x_\mu$.
Clearly, the elements $Z_{\fks\fkt}^\nu$ are linearly independent.

We now prove the spanning condition.
Suppose $h\in x_\la^\eta\sH'\cap\sH'x_\mu$. Then $h\in
 x_\la^\eta\sH'\pi_b\cap \pi_a^-\sH'x_\mu$, and,
by Cor. \ref{ind1}, 
$$h=\sum_{\nu,\bfs,\fkt}\al_{\bfs,\fkt}^\nu   
\pi_a^-X_{\bfs\fkt}^\nu=\sum_{\nu^\star,\fks,\bft} 
\be^{\nu^\star}_{\fks,\bft} \eta(X^{\nu^\star}_{\fks,\bft})\pi_b.$$
Equating the coefficients with respect to the basis
$\sY_{a^-,b}$, we see that  
$\al_{\bfs,\fkt}^\nu=\be^{\nu^\star}_{\fks,\bft}$
whenever $\bfs\in\bT_{\fks^\star}$ and
$\bft\in\bT_{\fkt^\star}$.
Consequently, $h$ is a linear combination of the elements $Z_{\fks\fkt}^\nu$'s.

To see the last assertion, we observe
from the third equality above that
$$\aligned
Z_{\fks\fkt}^\nu
&=\sum_{\bfs^\star\in\bT_{\fks^\star},\bft\in\bT_{\fkt}} 
T_{\bar\de(\bfs_2)^{-1}}T_{w_{a,b}\bar\de(\bfs_1)^{-1}w_{b,a}}
\bar x_{\nu^{(2)}}
\pi(a^-,b)\bar x_{\nu^{(1)}}T_{w_{b,a}\bar\de(\bft_2)w_{a,b}}T_{\bar\de(\bft_1)}\cr
&=\sum_{\bfs^\star\in\bT_{\fks^\star},\bft\in\bT_{\fkt}}
(T_{\bar\de(\bfs_2)^{-1}}\bar x_{\nu^{(2)}}T_{\bar\de(\bft_2)})
\pi(a^-,b)
(T_{\bar\de(\bfs_1)^{-1}}\bar x_{\nu^{(1)}}T_{\bar\de(\bft_1)})\cr
\endaligned
$$
which is  an element in
$\sY_{\la^-,\mu}$.
\end{pf}
 
\section{The homological property}

In this section, we will prove the homological property required
in stratifying the $q$-Schur$^{1.5}$ algebra.
The validity of such property is closely related to the existence of nice
bases as is seen in \cite{DPS1} where a Kazhdan-Lusztig basis was used
for {\it parabolic} subgroups.
We will see below that the bistandard basis and those they induce
for {\it quasi-parabolic} groups of type $B$ Hecke algebras
are sufficient to guarantee the homological property. 
We need some preparation. Recall from \cite[(1.2.1a)]{DS2}
the definition of the dominance order  $\trianglerighteq$ on bicompositions.
That is,
$\la\trianglelefteq\mu$ iff $\sum_{i=1}^j\laa_i\le\sum_{i=1}^j\mu^{(1)}_i$ for
all $j$ and $|\laa|+\sum_{i=1}^{j'}\lab_i\le
|\mu^{(1)}|+\sum_{i=1}^{j'}
\mu^{(2)}_i$
for all $j'$.

For $\mu\in\Pi^+_{2,b}$,
let 
$\Pi^+_\mu=\{\nu\in\Pi^+_2\mid \nu\trianglerighteq\mu\}$ and
$\Pi^+_{\mu,b}=\{\nu\in\Pi^+_{2,b}\mid \nu\trianglerighteq\mu\}$.
Note that $\Pi^+_\mu$ indexes the composition factors
of $\sH_Kx_\mu$, where $K$ is the quotient field of $\sZ$. 
We {\it linearly} order $\Pi^+_\mu$ by $\nu_{(1)}\preceq\nu_{(2)}\preceq\cdots$
such that $\nu_{(i)}\trianglerighteq\nu_{(j)}$ implies $i\le j$.
Fix an ordering on each $\fT^{ss}(\nu_{(i)},\mu)$ and concatenate
them via the linear ordering on the $\nu_{(i)}$ to obtain
a linear ordering $\fkt_1\preceq\fkt_2\preceq\cdots\preceq\fkt_{n_\mu}$ on 
 $\cup_{i=1}^{n_\la}\fT^{ss}(\nu_{(i)},\mu)$. 
Thus, using this ordering,
the basis $\sX_{\om_0,\mu}$  allows us to  define an integral {\it left} 
(twisted) ``Specht'' filtration 
for $\sH'\wx_\mu$ by setting: 
$$
\cases\sE_{i\sZ'}^\mu={\text{\rm span}}\{X_{\bfs\fkt_j}^{\nu(\fkt_j)}\mid 
1\le j\le i,
\bfs\in\bT^s(\nu(\fkt_j))\},\,\,\,
1\leq i\le n_\mu,\cr
\sE^\lambda_{0\sZ'}=0.\endcases$$
Here $\nu(\fkt_j)$ is the shape of $\fkt_j$, that is,
the partition whose Young diagram
is the underlying diagram of $\fkt_j$.
We shall write $\sE^\mu_\bullet$ for $\sE^\mu_{\bullet\sZ}$.

Note that the restriction of the ordering on $\Pi^+_\mu$ induces an 
ordering $\nu_{(i_1)}\preceq\nu_{(i_2)}\preceq\cdots$ on $\Pi^+_{\mu,b}$, 
and hence an ordering
 $\fkt_{j_1}\preceq\fkt_{j_2}\preceq\cdots\preceq\fkt_{m_\mu}$
on  $\cup_j\fT^{ss}(\nu_{(i_j)},\mu)$. The following lemma
is an easy consequence of Theorem \ref{ind2}.

\begin{lem} \label{tsf1}
Let $\la,\mu\in\Pi^+(r)$ such that $|\laa|+|\mu^{(1)}|=r$,
and let 
$$0=\sE_{0\sZ'}^\mu\subseteq
\sE^\mu_{1\sZ'}\subseteq\cdots\subseteq
\sE_{n_\mu\sZ'}^\mu=\sH'x_\mu$$
be the filtration defined as above. For each $i$, let
$j_{m(i)}$ be the maximal index in the sequence
$j_1,\cdots,j_{m_\mu}$ such that $j_{m(i)}\le i$.
Then
$x_\la^\eta\sH'\cap\sE^\mu_i$ is $\sZ'$-free for every $i$,
with basis
$$\sY_{\la^-,\mu, i}=\{Z_{\fks\fkt}^{\nu(\fkt)}\mid  
\fkt\in\{\fkt_{j_1},\cdots,\fkt_{j_{m(i)}}\},\fks^\star\in\fT^{ss}_b(\nu(\fkt),\la^\star)
\}.$$
\end{lem}

\begin{pf} Let $|\laa|=a$ and $|\mu^{(1)}|=b$.
By the definition above, $\sE^\mu_{i\sZ'}$ has basis
$$\{X_{\bfs\fkt_j}^{\nu(\fkt_j)}\mid 
1\le j\le i,\bfs\in\bT^s(\nu(\fkt_j))\}.$$
Thus, $\pi_a^-\sE^\mu_{i\sZ'}$ has basis
$$\{\pi_a^-X_{\bfs\fkt}^{\nu(\fkt)}\mid 
\fkt\in\{\fkt_{j_1},\cdots,\fkt_{j_{m(i)}}\},
\bfs^\star\in\bT^s_b(\nu(\fkt))\}.$$
Since $x_\la^\eta\sH'\cap\sE^\mu_{i\sZ'}=x_\la^\eta\sH'\pi_b\cap \pi_a^-
\sE^\mu_{i\sZ'}$,
it follows from a similar argument for Theorem
\ref{ind2} that $x_\la^\eta\sH'\cap\sE^\mu_{i\sZ'}$ is free 
with the required basis
$\sY_{\la^-,\mu, i}$.
\end{pf}

We will also need the following result which can be obtained by 
applying $\iota:T_w\mapsto T_{w^{-1}}$ to \cite[(6.1.6)]
{DS2}.

\begin{lem} \label{tsf2}
Keep the notation introduced above and
let
$\la,\mu\in\Pi^+_2$.
Then $x_\la\sH'\cap\sE_{i\sZ'}^\mu$ is free with basis
$$\sX_{\la,\mu, i}=\{X_{\fks\fkt_j}^{\nu(\fkt_j)}\mid 1\le j\le i,
\fks\in\fT^{ss}(\nu(\fkt_j),\la)\}.$$
\end{lem}

By taking duals, we turn the twisted Specht filtration $\sE^\mu_\bullet$
above
to a Specht filtration $\sF^\bullet_\mu$:
$$0=\sF_\mu^0\subseteq\sF_\mu^1\subseteq\cdots\subseteq\sF^{n_\mu}_\mu=x_\mu
\sH,$$
where
$$\sF^i_\mu=(\sH \wx_\mu/\sE^\mu_{n_\mu-i})^*,\,\, 0\leq i\leq n_\mu,$$ 
(compare \cite[(5.2.3)]{DS2}). We define the Specht module
$$S_\la=\sF^1_\la.$$

\begin{thm}\label{hp1}
For  $\la,\mu\in\Pi_2^+(r)$ and
any commutative $\sZ$-algebra
$\sZ'$, let $a=|\laa|$ and $b=|\mu^{(1)}|$. Then
base change defines isomorphisms
$$\cases
(1)\,\,
\Hom_{\sH}(\wx_\mu\sH/\sF^i_{\mu},\wx_\la\sH)_{\sZ'}
\overset{\sim}\to
\Hom_{\sH_{\sZ'}}((\wx_\mu\sH/\sF^i_{\mu})_{\sZ'},\wx_\la\sH_{\sZ'}),\cr
(2)\,\, \Hom_{\sH}(\wx_\mu\sH/\sF^i_{\mu},\wx_\la^\eta\sH)_{\sZ'}
\overset{\sim}\to
\Hom_{\sH_{\sZ'}}((\wx_\mu\sH/\sF^i_{\mu})_{\sZ'},\wx_\la^\eta\sH_{\sZ'}),\cr
\endcases
$$
for all $i$, assuming $a+b=r$ in the latter case. Also, the $\sZ'$-modules in \ref{hp1}(1) (resp. \ref{hp1}(2))
is free of rank $r_{\la,\mu, i}=\#\sX_{\la,\mu, i}$
(resp. $r_{\la^{\!-},\mu, i}=\#\sY_{\la^-,\mu, i}$).  
Furthermore, for $i<j$, the natural maps
$$\cases
(1')\,\,\Hom_{\sH}(\wx_\mu\sH/\sF^j_{\mu},\wx_\la\sH)_{\sZ'}
\to\Hom_{\sH}(\wx_\mu\sH/\sF_{\mu}^i,\wx_\la\sH)_{\sZ'},\cr
(2')\,\,\Hom_{\sH}(\wx_\mu\sH/\sF^j_{\mu},\wx_\la^\eta\sH)_{\sZ'}
\to\Hom_{\sH}(\wx_\mu\sH/\sF_{\mu}^i,\wx_\la^\eta\sH)_{\sZ'},\text{ if }a+b=r,\cr
\endcases
$$
have $\sZ'$-free cokernels of rank $r_{\la,\mu, i}-r_{\la,\mu, j}$
(resp. $r_{\la^{\!-},\mu, i}-r_{\la^{\!-},\mu, j}$).
 \end{thm}

\begin{pf}  Since
$(\wx_\mu\sH/\sF^i_{\mu})^*\cong\sE_{n_\mu-i}^{\mu}$,
$(\wx_\la\sH)^*\cong\sH \wx_\la$ (\cite[4.3]{DS2})
 and $(\wx_\la^\eta\sH)^*\cong
\sH \wx_\la^\eta$, (1) and (2) are equivalent to showing that base change
defines isomorphisms
$$\aligned
\Hom_{\sH}(\sH \wx_\la,\sE_i^{\mu})_{\sZ'}&
\overset{\sim}\to
\Hom_{\sH'}(\sH' \wx_\la,\sE_{i\sZ'}^{\mu}),\text{ and }\cr
\Hom_{\sH}(\sH \wx_\la^\eta,\sE_i^{\mu})_{\sZ'}&
\overset{\sim}\to
\Hom_{\sH'}(\sH' \wx_\la^\eta,\sE_{i\sZ'}^{\mu})\text{ if }a+b=r,\cr
\endaligned$$
for all $i$.
Since $\sH'$ is a Frobenius algebra over $\sZ'$ (and thus
``injective'' for $\sH'$-modules relative to $\sZ'$-split
exact sequences), we have 
$$\Hom_{\sH'}(\sH' \wx_\la,\sE_{i\sZ'}^{\mu})\cong x_\la\sH'\cap
\sE_{i\sZ'}^{\mu},$$ 
etc., reducing the isomorphisms to the isomorphisms
$$(x_\la\sH\cap\sE_{i}^{\mu})_{\sZ'}\overset{\sim}\to x_\la\sH'\cap\sE_{i\sZ'}^{\mu},
\text{ and } 
(x_\la^\eta\sH\cap\sE_{i}^{\mu})_{\sZ'}\overset{\sim}\to x_\la^\eta
\sH'\cap\sE_{i\sZ'}^{\mu}$$
if  $a+b=r$ in the latter case.
Now, Theorems \ref{tsf1} and \ref{tsf2} have established the above isomorphisms,
and the first assertion is proved.

The discussion in the previous paragraph shows that the 
maps in ($1'$) and ($2'$) 
are determined by the inclusions
$$x_\la\sH'\cap\sE_{n_\mu-j,\sZ'}^{\mu}\subseteq
x_\la\sH'\cap\sE_{n_\mu-i,\sZ'}^{\mu}
\text{ and }
x_\la^\eta\sH'\cap\sE_{n_\mu-j,\sZ'}^{\mu}\subseteq
x_\la^\eta\sH'\cap\sE_{n_\mu-i,\sZ'}^{\mu}.$$
Since the basis 
described in \ref{tsf1} or \ref{tsf2} for $n_\mu-j$ is a subset of 
the basis for $n_\mu-i$, the last assertion now follows easily..
\end{pf}

We observe that the isomorphism (2) in
(\ref{hp1}) holds for $a+b>r$ as 
both sides would be zero in this case. Also, one would expect that
the base change property in (\ref{hp1}(2)) holds
in general, though we do not need the case $a+b<r$.

We now can apply \cite[1.2.13]{DPS1} to obtain the following
homological vanishing property. Note that $\sH_K$, where
$K$ is the quotient filed of $\sZ$, is semisimple.

\begin{thm} \label{hp3} For $\la,\mu\in\Pi_2^+(r)$, we have
for all $i$,
$$\cases
(1)\quad
\Ext^1_{\sH}(\wx_\mu\sH/\sF^i_\mu,\wx_\la\sH)=0\cr
(2)\quad
\Ext^1_{\sH}(\wx_\mu\sH/\sF^i_\mu,\wx_\la^\eta\sH)=0\quad
\text{ if }|\laa|+|\mu^{(1)}|=r.
\cr
\endcases
$$
\smallskip
Moreover,  both $\Hom_{\sH}(\sF^j_{\mu},\wx_\la\sH)$
and $\Hom_{\sH}(\sF^j_\mu,\wx_\la^\eta\sH)$ (assuming 
$|\laa|+|\mu^{(1)}|=r$ in the latter case) are 
$\sZ'$-free for all $j$.
They have ranks $r_{\la,\mu}-r_{\la,\mu, j}$
and $r_{\la^{\!-},\mu}-r_{\la^{\!-},\mu, j}$, respectively, where
$r_{\la,\mu}=r_{\la,\mu, 0}=\#\sD_{\la,\mu}$ and
$r_{\la^{\!-},\mu}=r_{\la^{\!-},\mu, 0}=\#\sD_{\la,\mu}^0$.
\end{thm}
\begin{pf}
The first assertion 
follows from \cite[(1.2.13)]{DPS1} and the previous theorem. 

Next, observe that for any $j$, Theorem \ref{hp1}
implies that  $\Hom_{\sH}(\sF_{\mu}^j,
\wx_\la\sH)$ (resp.
$\Hom_{\sH}(\sF_{\mu}^j,
\wx_\la^\eta\sH)$) identifies with the cokernel of the natural map
$$\Hom_{\sH}(\wx_\mu\sH/\sF^j_{\mu},\wx_\la\sH)\to
\Hom_{\sH}(\wx_\mu\sH,\wx_\la\sH)$$
(resp. 
$\Hom_{\sH}(\wx_\mu\sH/\sF^j_{\mu},\wx_\la^\eta\sH)\to
\Hom_{\sH}(\wx_\mu\sH,\wx_\la^\eta\sH).$)
Thus, the last assertion follows from the last assertion of the previous 
theorem, together
with elementary properties of distinguished double coset representatives.
\end{pf}

\section{Stratifying the $q$-Schur$^{1.5}$ algebra}

We now turn back to the type $D$ case. So we assume in this section that
 $\sH'=\sH'_{1,q}$.
We have first the following base change property.

\begin{thm} \label{basec} Let $\sZ'$ be a commutative $\sZ$-algebra in which 2 is invertible.
Then the algebra $\sS_q^{1.5}(n,r;\sZ')$ is free and
 base change induces an isomorphism
$$\sS_q^{1.5}(n,r)_{\sZ'}\cong \sS_q^{1.5}(n,r;\sZ').$$
\end{thm}

\begin{pf} For $\la,\mu\in\Pi_{1.5}$, we have, using Lemma \ref{6b}c,
$$\aligned
     \Hom_{\dH'}(\dx_\mu\dH',\dx_\la\dH')
&\cong\Hom_{\sH'}(x_\mu\sH',\dx_\la\dH'\uparrow^{\sH'})\cr
&\cong\Hom_{\sH'}(x_\mu\sH',x_\la\sH')\oplus\Hom_{\sH'}(\dx_\mu\sH',x_\la^\eta\sH').
\cr\endaligned$$
Now, the Hom space $\Hom_{\sH'}(x_\mu\sH',x_\la\sH')$
has always the base change property by Theorem \ref{hp1}(1),
while the Hom space $\Hom_{\sH'}(\dx_\mu\sH',x_\la^\eta\sH')\cong x_\la^\eta\sH'\cap\sH'x_\mu
$
is  0 if  $|\laa|+|\mu^{(1)}|>r$  (Theorem \ref{piHpi3}), or
free of rank independent of $\sZ'$ if  $|\laa|+|\mu^{(1)}|=r$.
Hence the result follows.
\end{pf}

With this Theorem, we can determine the structure of the $q$-Schur$^{1.5}$
algebra immediately when $r$ is odd, at least relative to the $q$-Schur$^2$
algebra.

\begin{cor}\label{T:odd} Assume again that 2 is invertible in $\sZ'$ and $r$ is odd. Then $\sS_q^{1.5}(n,r;\sZ')$ 
is a centralizer subalgebra
of $\sS_q^2(n,r;\sZ')$ associated to the coideal $\Pi_{1.5}^+(n,r)$ of $\Pi_2^+(n,r)$. 
Hence, it is quasi-hereditary, and 
\begin{equation}
\sS_q^1(n,r;\sZ')\le\sS^{1.5}_q(n,r;\sZ')\le\sS^{2}_q(n,r;\sZ').
\end{equation}
\end{cor}

\begin{pf} When $r$ is odd, 
we have always $|\laa|+|\mu^{(1)}|>r$. It follows that 
$\Hom_{\sH'}(\dx_\mu\sH',x_\la^\eta\sH')=0$.
Therefore,
$$\Hom_{\sH'}(x_\la\sH',x_{\mu}\sH')\cong\Hom_{\dH'}(\dx_\la\dH',
\dx_{\mu}\dH'),$$
for all $\la,\mu\in\Pi_2$,
and consequently, $\sS^{1.5}_q(n,r;\sZ')$ is a centralizer
subalgebra of $\sS^{2}_q(n,r;\sZ')$ defined by
the coideal $\Pi_{1.5}^+$. The result follows. Note that
the left hand inequality has been proved in (\ref{clii}) above.
\end{pf}

When $r$ is even, 
it is very unlikely that $\sS_q^{1.5}(n,r;\sZ')$ is
quasi-hereditary in general.\footnote{A bad case is $r=4$, $q=1$, 
and $\sZ'$ the localization
$\BZ_{(2)}$. The five ordinary irreducible representations
associated to the four
bipartitions (2;2), (2;1,1), (1,1;2), (1,1;1,1) may be shown
to have only four distinct 2-modular constituents for
$\sS_q^{1.5}$. This precludes any integal
quasi-hereditary structure compatible with our partial order.}
However, it can be integrally and standardly stratified naturally, using
the type $B$ Specht modules filtrations $\sF_\la^\bullet$.
It suffices now, by the base change property,
to look at the integral case. Recall that $S_\la=\sF_\la^1$.

\begin{thm} \label{stra}
Let $\sZ'$ be the ring obtained by localizing $\sZ=\BZ[q,q^{-1}]$ at $2$.
For any $\la\in\Pi_{1.5}^+$, let $\De^{1.5}(\la)=
\Hom_{\dH'}(S_{\la\sZ'},\dT^{1.5}_{\sZ'})$.  Then $\De^{1.5}(\la)$ is $\sZ'$-free
and $\{\De^{1.5}(\la)\}_{\la\in\Pi_{1.5}^+}$ is a standard stratifying system 
for the category
of $\sS_q^{1.5}(n,r)_{\sZ'}$-modules.
\end{thm}

\begin{pf} The freeness follows from the definition
of $S_\la$ and Theorem \ref{hp3}
It remains to check the Hypothesis given in \cite[(1.2.9)]{DPS1}.

Let $\La=\Pi_{1.5}^+(n,r)$ and define 
a new order $\le$ on $\La$ by setting $\la\le\mu$ if
$\la\trianglelefteq\mu$ or $\la\trianglelefteq\mu^\star$ with
 $\la,\mu,\mu^\star\in\La$.
Clearly, $(\La,\le)$ is a quasi-poset. Let $\bar\La$ denote the 
associated poset.
 By \cite[(3.2.2c)]{DS2},
we have $\dT^{1.5}_{\sZ'}\cong\oplus_{\mu\in\La}\dT_\mu^{\prime\oplus d_\mu}$, where $\dT_\mu'=\dx_\mu\dH'\cong x_\mu\sH'$.

For $\mu\in\La$, $\dT_\mu'$ has a filtration $\sF_\mu^\bullet$ 
(see (\ref{hp1})) for which
$Gr^i \sF_\mu^\bullet\cong S_{\nu_{\mu,i}\sZ'}$ and $\nu_{\mu,i}>\mu$ if $i<n_\mu$. Note that this can be replaced
by $\bar \nu_{\mu,i}>\bar\mu$.

For $\la,\mu\in\La$, if $\Hom_{\dH'}(S_{\mu\sZ'},\dT_\la')\neq0$, then
$\Hom_{\dH_K}(\dT_{\la K},S_{\mu K})\neq0$. Thus, we have
either $\Hom_{\sH_K}(\sT_{\la K},S_{\mu K})\neq0$ or 
$\Hom_{\sH_K}(\sT_{\la K}^\eta,S_{\mu K})\neq0$ by Lemma \ref{6b}b. The former
implies that $\la\trianglelefteq\mu$, while the latter implies that
$|\laa|=|\lab|=|\mu^{(1)}|=|\mu^{(2)}|$ (thus, $\mu^\star\in\La$),
and $\la\trianglelefteq\mu^\star$. Therefore, $\la \le \mu$.

It remains to check the homological condition 
$$\Ext_{\dH'}^1(\dx_\mu\dH'/\sF_{\mu\sZ'}^i,\dx_\la\dH')=0$$
for all $i$ and $\la\in\Pi_{1.5}^+$.
However, one sees easily that 
$$\aligned
&\Ext_{\dH'}^1(\dx_\mu\dH'/\sF_{\mu\sZ'}^i,\dx_\la\dH')\cr
\cong&\Ext_{\sH'}^1(x_\mu\sH'/\sF_{\mu\sZ'}^i,x_\la\dH'\uparrow^{\sH'})\cr
\cong&\Ext^1_{\sH'}(\wx_\mu\sH'/\sF^i_\mu,\wx_\la\sH')\oplus
\Ext^1_{\sH'}(\wx_\mu\sH'/\sF^i_\mu,\wx_\la^\eta\sH')=0,
\cr\endaligned$$
by Theorem \ref{hp3}. Now, the theorem follows from
\cite[(1.2.10)(1.2.12)]{DPS1}.
\end{pf}

Let $\sO$ be regular local ring of Krull dimension $\le 2$ with the
field $K$ of fractions and residue field $k$. Assume further that
2 is invertible in $\sO$.

\begin{cor} \label{straa} The decomposition matrix of $\sS^{1.5}_q(n,r)_\sO$ 
contains an upper unitriangular block of size equal to the number of
irreducible modular representations of $\sS^{1.5}_q(n,r)_k$.
Moreover, if all $\De(\la)_k$ with $\laa=\lab$ are decomposable,
then $\sS^{1.5}_q(n,r)_k$ is quasi-hereditary.
\end{cor}

\begin{pf} Let $\widetilde A=\sS^{1.5}_q(n,r)$.
Since $\widetilde A_K$ is split semisimple,
we may assume that $\sO$ is complete.
By Theorem \ref{stra} and \cite[(1.2.5),(1.2.8)]{DPS1}, 
$A:=\widetilde A_k$ has a standard stratification
$$0=J_0\subset J_1\subset\cdots\subset J_m=A$$
such that (1) $m=\#\bar\La$; (2) $J_i/J_{i-1}=(A/J_{i-1})e_i(A/J_{i-1})$
for some idempotent $e_i$ with $(A/J_{i-1})e_i\cong\De(\la(i))_k$.
Moreover, there are $\widetilde A_\sO$-projective modules
$P(\la(i))$ such that $P(\la(i))$ has a $\De$-filtration (see \cite[(1.2.4)]
{DPS1}) with top section isomorphic to $\De(\la(i))_\sO$, 
other sections $\De(\mu)_\sO$ satisfying $\mu>\la(i)$, and such that
$P(\la(i))_k$ is the projective cover of $\De(\la(i))_k$.

Now, if $\la=\la(i)$ is not of the form
$(\al,\al)$, then $\De(\la(i))_k$ is indecomposable and 
$\End_A(\De(\la(i))_k)$ is a division ring, and hence
$J_i/J_{i-1}$ is a heredity ideal of $A/J_{i-1}$.
If $\la=(\al,\al)$ and   $S_{\la\sO}|_{\dH_\sO}$ 
is decomposable, then $S_{\la\sO}|_{\dH_\sO}$
a direct sum of two non-isomorphic indecomposable modules, 
and hence,  $\De(\la(i))_k$
is a direct sum of two distinct PIMs (of $A/J_{i-1}$).
Thus, in this case, $J_i/J_{i-1}$ is also a heredity ideal of $A/J_{i-1}$,
and the second assertion follows. 

Finally, if $\la=(\al,\al)$ and   
$S_{\la\sO}|_{\dH_\sO}$ is indecomposable,
then $\De(\la(i))_k$ is indecomposable and $\End_A(\De(\la(i))_k)$ is a local
ring (not a division ring!). Note that $\De(\la(i))_K$ is a sum of
two non-isomorphic simple modules.
In all cases,  we see that the matrix
$(d_{ij})$, where $d_{ij}=\dim\Hom_A(P(\la(i))_K,\De(\la(j))_K)$,
is upper unitriangular, and hence the result follows from \cite[(1.1.3)]{DPS1}.
\end{pf}

\begin{rems}\label{char2}  
(a) It is very likely that the odd characteristic condition
in the main results \ref{basec}-4 
can be removed. To do this we need a direct argument for type $D$
parallel to the work for type $B$ in \cite{DS2}. 
For example, when $r$ is odd, results (\ref{basec}) and (\ref{T:odd})
can be proved as follows:
Put
$$\dD_\la=\cases (\wsD_\la\cap \dW)\cup (s_0\wsD_\la\cap \dW),\text{ if }
|\laa|\ge1\cr
\wsD_\la \cap \dW,\text{ if }|\laa|=0.\cr\endcases
$$
Then one may prove that $\#\sD_{\la,\mu}=\#\dD_{\la\mu}$, where
$\la,\mu\in\Pi_{1.5}$ and 
$\dD_{\la\mu}=\dD_\la\cap\dD_\mu^{-1}$.
By directly constructing a basis for
$\Hom_{\dH'}(\dT_\mu',\dT_\la')$ (compare \cite[(4.2.6)]{DS2}), one
could prove that 
$$\Hom_{\sH'}(\sT_\mu',\sT_\la')=\Hom_{\dH'}(\sT_\mu',\sT_\la').$$
Thus, the quasi-heredity of the $q$-Schur$^{1.5}$ algebra
when $r$ is odd follows immediately.

(b) Note that (\ref{straa}) applies as well to
the $q$-Schur$^{2.5}$ and $q$-Schur$^{\oned}$ algebras in the 
linear prime case, because of (\ref{mrii}) above. Also, together
with \cite[7.15]{GrH}, we see that, in the linear prime case,
all $\De(\la)_k$ with $\laa=\lab$ are decomposable. Therefore,
by (\ref{straa}), the $q$-Schur$^{1.5}$ algebra is quasi-hereditary in this
case. It is likely that this result always holds when the polynomial
$\dg_r$ is invertible.
\end{rems}

\section{The bad prime case}

A case at ``the other extreme'' compared with the linear prime case
(\S3) is $p=2$, and hence, for finite groups of Lie type,  $q=1$ taking $q$ to be  prime power
with $p\nmid q$. (See the footnote below.) In this case, the group algebra
of $\dC=\langle t_1t_2,\cdots,t_1t_r\rangle$ has only one irreducible
modular representation, and is a local ring. (In the linear prime case
the Hecke algebra behaves as if it had a semisimple
subalgebra associated to $\dC$.) Nevertheless,
it is quite
easy to find the modular irreducible representation of
$\sS_q^{\oned}(r,r)$ (or $\sS_q^{\oneb}(r,r)$)
in terms of those for  $\sS_q^{1}(r,r)$ (which maybe regarded
as part of $\sS_q^{1.5}(r,r)$).
For simplicity, we will stick to the $\sS_q^{\oned}$ case.

To fix notation, we get $q_0=q=1$ and let $\sZ'$ be a commutative
local ring with residue field $k=\sZ'/\frak m$ of 
characteristic 2.\footnote{It is interesting to note that ``decomposition 
numbers'' for $\sS^{\oned}_{\sZ'}$, with $\sZ'=\BZ_{(2)}$
(or any char. 0 DVR with 2 in its maximal ideal), 
are the same for $q=1$ and for
$q$ equal to any odd prime power. This follows using, say,
\cite[(1.1.2)]{DPS1} and the fact that $\sS^{\oned}_{\BQ}$
is split semisimple for each of these specializations of $q$.}
Let $\sS^{\oned}_{\sZ'}$ denote
$\sS_q^{\oned}(r,r,\sZ')$ for some fixed $r>0$, and let 
$\sS_{\sZ'}^1$ denote $\sS_q^{1}(r,r,\sZ')$. The latter is just the
classical Schur algebra.

\begin{thm} \label{badp} There is a natural surjective homomorphism
$$\th:\sS_{\sZ'}^{\oned} \to \sS_{\sZ'}^1$$
whose kernel is contained in the radical $\rad(\sS_{\sZ'}^1)$.
In particular, the irreducible representations of $\sS_{\sZ'}^{\oned}$
over $k$ are all obtained from those $\sS_{k}^1$
by factorization through $\th$.

Also, if $\sZ'\to\sZ''$ is a commutative local rings, we have the base 
change properties
$$\sZ''\otimes \sS_{\sZ'}^{\oned}\cong \sS_{\sZ''}^{\oned}
\text{ and }
\sZ''\otimes \sS_{\sZ'}^1\cong \sS_{\sZ''}^1.
$$
\end{thm}

\begin{pf} Write $\dW=\dC\rtimes\oW$ where $\dC$ is the normal
elementary abelian 2-group $\langle t_1t_2,\cdots,t_1t_r\rangle$.
If $\oW_\mu$ is any subgrup of $\oW$, and $M$ is any right
$\sZ'\oW_\mu$-module, then the semidirect product decomposition gives
$$\Ind_{\sZ'\oW_\mu}^{\sZ'\dW}(M)\cong
\Ind_{\sZ'\oW_\mu}^{\sZ'\oW}(M)\otimes \sZ'\dC$$
where the right hand side of the tensor product is regarded
as a $\sZ'\dW$-module by inflation, and $\sZ'\dC$ is regarded
as a $\dW$-module with $\dC$ acting by right multiplication
and $\oW$ acting by conjugation. Thus,
$$\dT_{\sZ'}^{\oned}\cong  \bar\sT_{\sZ'}\otimes\sZ'\dC$$
with $\dT_{\sZ'}^{\oned}$ as in (\ref{typeD}),
and $\bar\sT_{\sZ'}$ defined similarly for the classical symmetric
group $\oW$. Passing to the endomorphism algebras, we have
$$
\sS_{\sZ'}^{\oned}=\End_{\sZ'\dW}(\dT_{\sZ'}^{\oned})
\cong\left(\End_{\sZ'}(\bar\sT_{\sZ'})\otimes \sZ'\dC\right)^\oW
$$
where the right hand side is the fixed points for an action of $\oW$ on the 
tensor product. Note that the augmentation
$\sZ'\dC\to\sZ'$ is a $\sZ'\oW$-homomorphism, split by the
natural inclusion
$\sZ'\cong\sZ\cdot 1\subseteq \sZ'\dC$.
If $I(\sZ'\dC)$ denotes the augmentation ideal of
$\sZ'\dC$, then a power of $I(\sZ'\dC)$ is contained in $2\sZ'\dC$.
The same is true for the ideal
$$\left(\End_{\sZ'}(\bar\sT_{\sZ'})\otimes I(\sZ'\dC)\right)^\oW$$
of
$\left(\End_{\sZ'}(\bar\sT_{\sZ'})\otimes\sZ'\dC\right)^\oW
\cong\sS_{\sZ'}^{\oned}$.
The quotient of $\sS_{\sZ'}^{\oned}$ by the ideal
is
$$\left(\End_{\sZ'}(\bar\sT_{\sZ'})\otimes\sZ'\right)^\oW\cong\sS_{\sZ''}^1.$$
This gives the surgective homomorphism
$\th$ described in the theorem.

The base change properties follows easily from the fact that both
$\sT_{\sZ'}^{\oned}$ and $\bar\sT_{\sZ'}$ are permutation modules
for $\sZ'\oW$. The proof is complete.
\end{pf}

\begin{rems}
(a) The $q$-analogue of the homomorphism $\th$ always exists.
This is because $\bar\sT_{\sZ'}=\pi_r^-\sT^{\oneb}_{\sZ'}$ (or, in the
type $D$ case, $\bar\sT_{\sZ'}=\dpi_r\sT^{\oned}_{\sZ'}$) and $\th$
is simply the restriction map. 
In the type $B$ case, we even know that it is onto and
what the kernel of $\th$ is.
Using the bistandard (or celluar) basis for $\sS^{\oneb}_q(n,r;\sZ')$
(see \cite[(6.1.1)]{DS2}), we see easily that
the basis elements $\Phi_{\fks\fkt}^\la$ with $|\laa|=0$  are sent to a
basis for $\sS_q^1(n,r;\sZ')$, while those with $|\laa|>0$ are sent to 0
under $\th$. The surjectivity in the type $D$ case
should follow from a direct construction of a standard basis for
$\sS^{\oned}_q$,
cf. (\ref{char2}).

(b) The bad prime case for finite groups of Lie type generalizes, at the
Hecke algebra level, to the case where 
the factor $q_0+1$ of $g_r$ (or the factor 2 of $\dg_r$ for type $D$)
is zero in the base ring. However, if $q\neq1$, the kernel of $\th$ defined
in (b) is not nilpotent in general. For example, in the type $B_2$ case,
consider the cellular basis elements at the $\la=(1;1)$ level. They are
$1+T_{s_0},$ $(1+T_{s_0})T_{s_1}$, $T_{s_1}(1+T_{s_0})$ and
$T_{s_1}(1+T_{s_0})T_{s_1}$. Since
$$\aligned
(1+T_{s_0})T_{s_1}(1+T_{s_0})T_{s_1}
&=(1+T_{s_0})T_{s_1}^2+(1+T_{s_0})T_{t_2}\cr
&=(1+T_{s_0})T_{s_1}^2-q(1+T_{s_0})+\pi_2\cr
&=(q-1)(1+T_{s_0})T_{s_1}+\pi_2,\cr
\endaligned
$$
We see that the Hecke algebra $\sH'$ has three distinct irreducible modular
representations, while $\bar\sH'$ has only two. It is easy to check that the same is true for
$\sS_q^{\oneb}(2,2;\sZ')$ and $\sS_q^{1}(2,2;\sZ')$.

 (c)
It is resonable to speculate, based on the section 3 and 8, that the irreducible
modular representations of $\sS_q^{\oned}$, in any characteristic
$p$, are determined by some of the irreducible
modular representations of 
$\sS_q^{1.5}$, at least when $q$ itself is a power of a different prime.
In the linear prime case, all the irreducible modular representations
are required, while only those associated to
$\sS_q^1\le\sS_q^{1.5}$ are needed for $p=2$ and $q=1$.
It would be very satisfactory if one could
predict just from the order of $q$ modulo $p$ just what part
of $\sS_q^{1.5}$ were required. We hope to explore possible
general theories along these lines in the future.

\end{rems}


\begin{thebibliography}{9}
\bibitem{C} M. Cabanes, {\em Alg\`ebres de Hecke comme alg\`ebres
symm\'etriques et th\'eor\`em de Dipper.},
C.R. Acad. Sci. Paris, {\bf 327} (1998), 531-536.

\bibitem{CPS}
E. Cline, B. Parshall and L. Scott, {\em
Stratifying endomorphism algebras}, Memoirs Amer. Math. Soc.
{\bf 591}, 1996.

\bibitem{DJ}
 R. Dipper and G. James, {\em Representations of Hecke algebras of type $B_n$}
J. Algebra, {\bf 146} (1992), 454--481.

\bibitem{DJM}
 R. Dipper, G. James and A. Mathas, {\em The $(Q,q)$-Schur algebra},
Proc. London Math. Soc. {\bf 77} (1998), 327--361.

\bibitem{Du} J. Du, {\em Cells in certain sets of matrices}, 
T\^ohoku Math. J. {\bf 48} (1996), 417--427.

\bibitem{Du1} J. Du, {\em Generalized $q$-Schur algebras and the ways
to approach them}, Virginia Conference Proceedings (to appear).

\bibitem{DPS1} J. Du, B. Parshall and L. Scott, {\em 
Stratifying endomorphism algebras
associated to Hecke algebras}, J. Algebra, {\bf 203} (1998), 169--210.

\bibitem{DPS3} J. Du, B. Parshall and L. Scott, {\em 
Quantum Weyl reciprocity and tilting modules}, Commun. Math. Phys.
{\bf 195} (1998), 321--352.

\bibitem{DS2} J. Du and L. Scott, {\em The $q$-Schur$^2$ algebra}, 
Trans. Amer. Math. Soc. (to appear).

\bibitem{GH} M. Geck and G. Hiss, {\em 
Modular representations of finite groups of 
Lie type in non-defining characteristics}, in {\it Finite reductive groups},
M. Cabanes, ed. (1996) 173--227.

\bibitem{GrH} J. Gruber and G. Hiss, {\em 
Decomposition numbers of finite classical
groups for linear primes}, J. reine angew. Math. 485 (1997), 55-91.

\bibitem{HW}  J. Hu and J.-p. Wang, {\em Hecke algebras of type $D_n$ at roots of unity}, J. ALgebra 212 (1999), 132-160.

\bibitem{Ma} C. Mak, {\em Quasi-parabolic subgroups of $({\Bbb Z}/ m {\Bbb Z})\wr
{\frak S}_r$}, preprint, UNSW.

\bibitem{M}  G. Murphy, {\em The representations of Hecke algebras
of type $A_n$}, J. Algebra {\bf 173} (1995), 97-121.

\bibitem{P}  C. Pallikaros, {\em Representations of Hecke algebras
of type $D_n$}, J. Algebra {\bf 169} (1994), 20-48.


\end{thebibliography}
\end{document}